\documentclass[10pt,reqno]{amsart}
\usepackage{amsmath}
\usepackage{amssymb,amscd}
\usepackage[all]{xy}
\usepackage[english]{babel}

\usepackage{times}
\usepackage{t1enc}
\usepackage{graphicx}
\usepackage{float}
\usepackage{tikz}
\usepackage{tikz-cd}
\usepackage{caption}
\usepackage{subcaption}
\usepackage{epstopdf}
\usepackage{epsfig}
\usepackage{mathrsfs}
\usepackage{setspace} 
\usepackage[mathscr]{euscript}
\usepackage[utf8]{inputenc}
\usepackage[pdftex]{hyperref} 
\usepackage{amsmath, amssymb, amsthm, amsfonts, amstext,  textcomp, bm, amscd}
\usepackage[all]{xy}
\usepackage[english]{babel}
\usepackage{xcolor,soul}
\usetikzlibrary{matrix,arrows}

\newtheorem{theorem}{Theorem}[section]
\theoremstyle{plain} 
\theoremstyle{plain} 
\newcommand{\thistheoremname}{}
\newtheorem{genericthm}[theorem]{\thistheforemname}

\newtheorem*{genericthm*}{\thistheoremname}
\newenvironment{namedthm*}[1]
{\renewcommand{\thistheoremname}{#1}%
	\begin{genericthm*}}
	{\end{genericthm*}}

\numberwithin{equation}{section}
\newtheorem{thm}{Theorem}[section]
\newtheorem*{thm*}{Theorem}
\newtheorem{prop}[thm]{Proposition}
\newtheorem{lem}[thm]{Lemma}

\newtheorem{rmk}{Remark}

\numberwithin{equation}{section}

\DeclareMathOperator{\Gal}{Gal}
\def\NN{\widetilde {N}_{\mathcal{O}_{H}}(4,0)}
\def\E{\mathbb{E}}
\def\R{\mathfrak{R}}

\def\F{\mathbb{F}}
\def\Z{\mathbb{Z}}
\def\A{\mathcal{A}}

\def\H{\mathcal{H}}
\def\C{\mathbb{C}}

\def\u{N_{q}(\text{M}^{s})}

\def\h{\H_{d,q}}
\def\bigcirc{O}
\def\w{N_q(M_{\text{L}}(n,d))}
\def\cM{\mathcal{M}_{L}(n,d)}
\def\M{M_{L}(n,d)}
\def\Mus{\mathcal{M}_{L}^{\text{us}}(n,d)}
\def\A{\frac{1}{N_{q}(\Aut E)}}
\def\fP{\mathfrak P}
\def \N{\bf \widetilde{N}}

\def\hh{\H_{\gamma,q}}

\def\a{\text{Higgs}}
\def\b {\text{Higgs}_{2,d}(X)}
\def\vv{	\log N_{q} \left(\text{Higgs}_{2,d}(X)\right)}
\def\vw{	\log N_{q} \left(\text{Higgs}_{2,d}(H)\right)}

\def \({\left(}
\def \){\right)}
\def \<{\langle}
\def \>{\rangle}
\def \bar{\overline}
\def \deg{\mathrm{deg}}

\def \Aut{{\rm Aut}}
\def \End{{\rm End}}

\def \Hom{{\rm Hom}}

\begin{document}
	
	\title[Statistics of Moduli Space]{Statistics of Moduli Spaces of vector bundles over hyperelliptic curves}

	\author[A. Dey]{Arijit Dey}
	\address{Department of Mathematics, Indian Institute of Technology-Madras, Chennai, India }
	
	\email{arijitdey@gmail.com}
	
	\author[S. Dey]{Sampa Dey}
	\address{Stat Math Unit,Indian Statistical Institute, Kolkata, India }
	
	\email{sampa.math@gmail.com}

	\author[A. Mukhopadhyay]{Anirban Mukhopadhyay} 
	\address{Institute of Mathematical Sciences, HBNI, CIT Campus, Taramani,  Chennai, India}
	
	\email{anirban@imsc.res.in}
	
	

	\begin{abstract}

		We give an asymptotic formula for the number of $\F_{q}$-rational points over a fixed determinant moduli space of stable vector bundles of rank $r$ and degree $d$ over a smooth, projective curve $X$ of genus $g \geq 2$ defined over $\F_{q}.$
		Further, we study the distribution of the error term when $X$ varies over a family of hyperelliptic curves. We then extend the results to the Seshadri desingularisation of the moduli space of semi-stable vector bundles of rank $2$ with trivial determinant, and also to the moduli space of rank $2$ stable Higgs bundles.
\end{abstract}

\maketitle

\tableofcontents

\section{Introduction}
Let $V$ be a quasi-projective variety defined over a finite field $\F_{q},$ and $\bar{V}:=V \otimes _{\F_{q}} \bar{\F_{q}}.$ We denote the cardinality of the set of $\F_{q}$-rational points on $\bar{V}$ by $N_{q}(V).$ Studying the quantitative behaviour of $N_q(V)$ is of paramount importance across several mathematical domains, such as finite field theory, number theory, algebraic geometry and so on.

Let $X$ be a smooth, projective curve of genus $g \geq 2$ over a finite field $\F_{q}$ such that $\bar{X}=X \otimes_{\F_{q}} \bar{\F_{q}}$ is irreducible, and $L$ be a line bundle on $X$ of degree $d,$ defined over $\mathbb F_q.$ Let $M(r,d)$ (resp, $M^s(r,d)$) be the moduli space of semistable (resp, stable) vector bundles of rank $r$ and degree $d,$ and $M_{L}(r,d)$ (resp. $M^{s}_{L}(r,d)$) be the moduli space of semistable (resp. stable) vector bundles of rank $r$ with fixed determinant $L.$ When $r$ and $d$ are coprime, the moduli space $M_{L}(r,d)$ is an irreducible smooth projective variety of dimension $(r^{2}-1)(g-1).$ Replacing $\F_q$ by a finite extension if necessary, we may assume that everything i.e. $X$, $L$ and $M_{L}(r,d)$ are defined over $\F_q.$ It is known that when $gcd(r,d)=1,$ the $\F_{q}$-rational points of $M_{L}(r,d)$ are precisely the isomorphism classes of stable vector bundles on $X$ defined over $\F_{q}$  [Proposition 1.2.1, \cite{HaNa}].
Now, when rank $r=1$ and degree $d=0,$ the moduli space $M(1,0)$ is the Jacobian $J_{X}$ of the curve $X,$ which is an abelian variety of dimension $g.$ Due to \textit{Weil conjectures}, the functional equation and analogue of the Riemann hypothesis for zeta function of a smooth projective curve of genus $g$ implies:
\[(\sqrt{q}-1)^{2g}\leq N_{q}(J_{X})\leq (\sqrt{q}+1)^{2g}.\]
For $g=1,$ this bound is tight due to the classical result of Deuring \cite{Deu}. For higher genus we know several improvements of this bound by Rosenbloom and Tsfasman\cite{RoTs}, Quebbemann\cite{Qu},  Stein and Teske\cite{StTe} and others. In \cite{Ts}, Tsfasman has shown that for a fixed finite field $\F_{q},$
\[g\log{q} + o(g)\leq \log(N_{q}(J_{X}))\leq g\left( \log{q}+(\sqrt{q}-1)\log\frac{q}{q-1}\right)+o(g) \]
as $g\rightarrow \infty.$ In terms of \textit{gonality} (that is the smallest integer $d$ such that $X$ admits a non-constant map of degree $d$ to the projective line over $\F_{q}$), 
Shparlinski \cite{Sh}, showed that  
\[ \log(N_{q}(J_{X}))=g\left(\log{q}+O_{q}\left(\frac{1}{\log{\left( \frac{g}{d}\right) }} \right)  \right) \]
as $q$ is fixed and $g\rightarrow \infty.$
In \cite{XiZa}, when the function field $\F_{q}(X)$ is a geometric Galois extension over the field of rational functions $\F_{q}(x)$ of degree $N,$ Xiong and Zaharescu estimated $N_{q}(J_{X})$ in terms of $q,\,g$ and $N.$ They gave the following explicit bound [Theorem 1,  \cite{XiZa}]
\[
| \log(N_{q}(J_{X}))-g\log{q} | 
\leq
(N-1)
\left( 
\log \text{max} 
\left\lbrace 
1, \frac{       \log{    \(    \frac{7g}{N-1}      \)    }     }             {\log q}
\right\rbrace 
+3 
\right) ,
\]
which holds true for any $q$ and $g.$ More precisely, we see the quantity $\left( \log(N_{q}(J_{X}))-g\log{q}\right) $ is essentially bounded by $O\left(  \log \log g  \right), $ which is significantly smaller than the bound $O\left(  \frac{g}{\log g} \right) $ given by Shparlinski in \cite{Sh}.

Motivated by the work of Xiong and Zaharescu \cite{XiZa}, we got interested in studying similar problems for moduli space of  stable rank $r$ and degree $d$ vector bundles with a fixed determinant  on a smooth projective curve. This can be interpreted as a non-abelian analogue of the work by Xiong and Zaharescu \cite{XiZa}. Now onwards, we will be considering a smooth projective curve $X$ of genus $g,$ defined over $\F_{q},$ where the function field $\F_{q}(X)$ is a geometric Galois extension over the field of rational functions $\F_{q}(x)$ of degree $N.$ We will simply call such curves as \emph{Galois curve} of degree $N.$ Here ``geometric '' means $\F_{q}$ is algebraically closed inside $\F_{q}(X).$
In a previous work, we explicitly studied the case when rank $r=2$ and degree $d=1,$ and we write the bound for $\log \left( N_{q}\left( M_{L}(2,1)\right) \right) $ in terms of $q,\,g$ and $N$ [Theorem 1.1, \cite{ASA}]. To prove this we used the Siegel formula \eqref{Si2.1}, where we need to count the number of isomorphism classes of unstable bundles too.  In this paper, we have generalized our previous results for the fixed determinant moduli space of any rank $r$ and degree $d$ with the condition that $gcd(r,d)=1.$ The challenging part was estimating the number of $\mathbb F_q$-rational points of  isomophism classes of unstable vector bundles, since the automorphism group is not constant. In this regard, the significant input comes from the work of Desale and Ramanan [Proposition 1.7, \cite{DeRa}], by which we could give an asymptotic bound in terms of $g,$ and $q$ inductively on the rank $r$ (cf. Proposition \ref{Pdom5.1}). It is worth mentioning here that one can do similar computations for the case of non-coprime $r$ and $d,$ too but the difficulty will arise in computing the number of $\mathbb F_q$-rational points 
of strictly-semistable strata.  In the rank $2$ and trivial determinant this was done in \cite{BaSe1} with great bit of care, and by using this we are able to do similar study for $M_{\mathcal O_X}(2,0)$ and it's Seshadri desingularization. 
Our first main result is the following asymptotic formula for $\w$ in terms of $N, \, q$ and $g.$ 

\begin{thm}\label{thm1.3}
Let	$X$ be a Galois curve of degree $N$ of genus $g\geq 2$ over $\mathbb{F}_{q}.$ Assume that $r$ and $d$ are coprime. 
If $\log{g} > \kappa\log(Nq)$ for some $\kappa>0$, sufficiently large absolute constant (independent of $N$), then for a constant $\sigma>0,$ depending only on the rank $r$ of the vector bundle 
we have

	\begin{align*}
	\log(\w)&= (r^2-1)(g-1)\log{q}\,
	 + 
O\left( 
A+q^{-\sigma g}\exp(A)
\right), 
	\end{align*}
	where 
	$A=N\left( \frac{1}{\sqrt q}+\frac{\log\log g}{q^2}\right) ,$ and 
	the implied constant depends on $r$ and $N.$
\end{thm}

Next we restrict our attention to the family of hyperelliptic curves i.e when $N=2$. Assume that $q$ is odd and  $\gamma$ is a positive integer $\geq\,5$.  Let $\H_{\gamma,q}$ be a family of curves given by the equation $y^2=F(x)$, where $F$ is a monic, square-free polynomial of degree $\gamma$ with coefficients in  $\mathbb{F}_{q}.$ Every such curve corresponds to an affine model of a unique projective hyperelliptic curve $H$, with genus $ g = \left[\frac{\gamma-1}{2}\right]$.
On $\H_{\gamma,q}$ we consider the uniform probability measure. With this set up, when $g$ is 
fixed and $q$ is growing, Katz and Sarnak showed that $ \sqrt{q}(\log N_{q}(J_{H})- g \log {q})$ is distributed as the 
trace of a random $2g\times2g$ unitary symplectic matrix \cite[Chapter 10, Variant 10.1.18]{KaSa}. On the other side, 
when the finite field is fixed and the genus $g$ grows, Xiong and Zaharescu \cite{XiZa} found the limiting distribution 
of $\log N_{q}(J_{H})- g \log {q}$ 
in terms of its characteristic function. Moreover, when both $g$ and $q$ grow, they showed that $ \sqrt{q}(\log N_{q}(J_{H})- g \log {q}) $ has a standard Gaussian distribution \cite{XiZa}.  Now, for every $H$ in $\H_{\gamma,q}$, fix a polarization i.e. a line bundle $L_{H}$ of degree $1$. In our previous work, we have studied the distributions of the quantity
$\log N_q (M_{L_H}(2,1)) - 3(g-1)\log q$ as the polarized curve $(H, L_H)$ varies over a large family of polarized hyperelliptic
curves. Since we are interested in $\mathbb F_{q}$-rational points of moduli spaces which is independent of the determinant  [Proposition 1.7, \cite{DeRa}], we denote the family of polarized curves by the same notation $\H_{\gamma,q}$. In this case, first we write the limiting distribution of $N_q (M_{L_H}(2,1)) -  3(g-1)\log{q}$ as $g$ grows and $q$ is fixed, in terms of it's characteristic function. Further, when $g$ and $q$ both grows together we see that $q^{3/2}\left(\log N_q (M_{L_H}(2,1)) -3(g-1)\log{q}\right)$ has a standard Gaussian distribution (see [Theorem 1.2, \cite{ASA}].) \\

When, $r\geq 2,$ and for any degree $d$, we consider the random variable, 
\[
\mathcal R_{(r,d)} : \H_{\gamma,q} \rightarrow \mathbb{R}
\]
given by, 
\begin{equation}\label{def:random variable}
\mathcal R_{(r,d)}(H):=
 \log N_{q}(M_{L_{H}}(r,d))- (r^2-1)(g-1)\log{q}.
\end{equation}
Before stating our next result here we fix some relevant notations. Over the polynomial ring $\F_{q}[t],$ we denote the degree of a polynomial $f$ by $\deg(f),$ and define the norm of a polynomial $|f|$ as $q^{\deg (f)}.$ Also, for any integer $k,$ let $\delta_{k/2}=1$ if $k$ is even and $0$ otherwise. Moreover, for a fixed $q,$ and a fixed rank $r\geq 2,$ we define
\begin{equation}\label{def:log constant}
C_{q}(r):=\log{\left(\frac{q^{(r^2-1)}}{\prod\limits_{k=2}\limits^{r}(q^{k-1}-1)(q^k-1)}\right)}
-\delta_{\gamma/2}\sum\limits_{k=2}\limits^{r}\log(1- 1/q^{k}).
\end{equation}

\begin{thm}\label{thm1.4}
	Let $(r, d) \in \mathbb{N} \times \mathbb{Z},$ such that $r\geq 2,$ and $gcd(r,d)=1.$ Then there exists an absolute constant $c>0$ such that
	\begin{equation*}
		\mathcal{R}_{(r, d)}-C_{q}(r)
		=
		\sum\limits_{k=1}^{r-1} \R_{(\gamma, q)}^{(k)}
		+ O(q^{-cg}),
		\end{equation*}
		  where $\R_{(\gamma, q)}^{(k)}$'s are random variables on $\hh$ satisfying the following:
		  \begin{itemize}
		  \item[(1)] For a fixed $q$ and $g=\left[\frac{\gamma-1}{2} \right] \rightarrow\infty,$ the random variable $\R_{(\gamma, q)}^{(k)}$ converges weakly to $\R_{k},$
		 where the characteristic function 
		$\phi_{\R_k}(t)=\mathbb{E}(e^{it\R_{k}}) $ of $\R_{k}$ is given by 
		\begin{align*}
			\phi_{\R_k}(t)
			&=
			1+ 
			\sum\limits_{n=1}\limits^{\infty}\frac{1}{2^{n}n!}
			\sum\limits_{\substack{P_{1},..., P_{n}}}\prod\limits_{j=1}\limits^{n}
			\left(\frac{(1-\mid P_{j} \mid^{-(k+1)})^{-it}+ (1+\mid P_{j} \mid^{-(k+1)})^{-it}-2}{(1+ \mid P_{j}\mid ^{-1})}\right),  \, \, 
		\end{align*}
		for any real number $t.$
		The inner sum is over distinct monic, irreducible polynomials $P_{1},...,P_{n}$ in $\F_{q}[x].$ 
		
		\item [(2)] 
		For a fixed $q,$ we have
	 \[
	 \lim\limits_{\gamma \rightarrow \infty}
	 Cov
	 \left( 
	 \R_{(\gamma, q)}^{(i)},\, \R_{(\gamma, q)}^{(j)}
	 \right) =\frac{1}{q^{i+j+1}} + O\left( \frac{1}{q^{i+j+2}}\right),
	\]
for any $1\leq i \neq j \leq r-1,$ 
	 
	 \item[(3)] If both $q, g\rightarrow\infty,$ then 
	${q^{\frac{(2k+1)}{2}}\R_{k}}$ has a Gaussian distribution. 

\end{itemize}
\end{thm}

\begin{rmk}\label{rmk1.2}
	When $r=2,$ from Theorem\ref{thm1.4}, we recover \cite[Theorem $1.2$]{ASA}.
\end{rmk}

\begin{rmk}\label{rmk on quotient}
	Let $T_{r}$ be the group of $r$-torsion points of the jacobian $J_{X}.$ Suppose the characteristic of the field $F_{q}$ is coprime with $r.$ 
Then over $\F_{q},$ the group scheme $T_{r}$ acts on the moduli space $M_{L}(r,d)$ by tensorisation.
Now $T_{r}$ being finite, the quotient $M_{L}(r,d)/T_{r}$ exists, and is a connected component of projective $PGL(r)$-bundles and is a projective variety. Further from \cite[Theorem 2]{HaNa}, it follows that
	\[
	\w=N_{q}\bigg( M_{L}(r,d)/T_{r} \bigg).
	\] 
	Using the above identity, Theorem \ref{thm1.3} and Theorem \ref{thm1.4} can be proved mutatis-mutandis for the quotient $M_{L}(r,d)/T_{r}$ which can be considered 
	as one of the component of moduli of $PGL(r)$-bundles.
	
\end{rmk}

Now, if we look at the case when rank and degree of the vector bundle are not coprime, the moduli space may not  be smooth. In particular, if we see the moduli space $M_{\mathcal O_X}(2,0),$ this is smooth only when genus of the curve $X$ is two [cf. \cite{NaRa1}]. For genus $\geq 3$, Seshadri constructed a natural desingularisation (moduli theoretic) of $M(2,0)$ over any algebraically closed field $k,$ having characteristic other than $2$ (cf.\cite{Se}).
He constructed a smooth projective variety $N(4,0),$ whose closed points corresponds to the $S$-equivalence classes of parabolic stable vector bundles of quasi-parabolic type (4,3) together with small weights $(\alpha_1,\alpha_2)$ such that the underlying bundles are semi-stable of rank $4$ and degree $0$ and having the endomorphism algebras as specialisations of $(2\times 2)$- matrix algebras.The natural desingularisation map from $N(4,0)$ to $M(2,0)$ is an isomorphism over $M^{s}(2,0)$ (see [Theorem 2, \cite{Se}] for more details). 
Furthermore, restricting on the subvariety of $N(4,0)$ whose closed points corresponds to isomorphism classes of parabolic stable vector bundles with the determinant of underlying bundle isomorphic to $\mathcal{O}_{X},$ is the desingularisation of $M_{\mathcal{O}_{X}}(2,0)$ (cf. Theorem $2.1$ in \cite{BaSe1}), and we denote this moduli space by ${\widetilde {N}}_{\mathcal{O}_{X}}(4,0)$.

In this paper, over the family of hyperelliptic curves $\H_{\gamma,q}$, we study the distribution of the $\F_{q}$-rational points on $M^{s}_{\mathcal{O}_{H}}(2,0)$ and its Seshadri disingularisation model ${\NN}.$ First we define the random variable
\[
\mathcal R_{(2,0)} : \H_{\gamma,q} \rightarrow \mathbb{R}
\]
such that
\begin{equation*}\label{def:random variable1}
\mathcal R_{(2,0)}(H):=
\log N_{q}(M_{\mathcal{O}_{H}}^{s}(2,0))- 3(g-1)\log{q}.
\end{equation*}
We have following similar results for $M^{s}_{\mathcal{O}_{H}}(2,0).$

\begin{thm}\label{thm1.5}
	\begin{itemize}
		\item[(1)]
		If $q$ is fixed and $g\rightarrow\infty,$ then 
	\[
	\mathcal{R}_{(2,0)}-C_{q}(2)
	\]
converges weakly to a random variable $\R,$ whose characteristic function $\phi_{\R}(t)$ is given by 
	\begin{eqnarray*}
		\phi_{\R}(t)&=&
		 1+ 
		 \sum\limits_{n=1}\limits^{\infty}
		 \frac{1}{2^{n}n!}\sum\limits_{\substack{P_{1},...,P_{n}}}
		 \prod\limits_{j=1}\limits^{n}\left(\frac{(1-\mid P_{j} \mid^{-2})^{-it}
		 	+ (1+\mid P_{j} \mid^{-2})^{-it}-2}{(1+ \mid P_{j}\mid ^{-1})}\right),  
	\end{eqnarray*}
	for all real number $t,$ where the inner sum is over distinct monic irreducible polynomial $P_{1},...,P_{n}$ in $\F_{q}[t].$
	
	\item[(2)]
	 If both $q, g\rightarrow\infty,$ then $q^{3/2}\mathcal{R}_{(2,0)}$ has a standard Gaussian distribution.  

\end{itemize}
\end{thm}

\begin{rmk}\label{rmk1.3}
	For the moduli space $M_{\mathcal{O}_{H}}(2,0),$ 
	when $q$ is fixed and $g\rightarrow\infty,$ over the family $\hh$ one can see that the random variable
	\[\log N_{q} (M_{\mathcal{O}_{H}}(2,0))  -3(g-1)\log{q}-C_{q}(2)
	\] 
	converges weakly to a random variable $\R,$ whose characteristic function is the same 
	$\phi_{\R}(t)$
	as given in Theorem \ref{thm1.5}. 
	Furthermore, when both $q, g\rightarrow\infty,$ then the random variable
	\[
	q^{3/2}
	\left(
	\log N_{q}( M_{\mathcal{O}_{H}}(2,0)) -3(g-1)\log{q}
	\right)
	\] 
	has a standard Gaussian distribution. 
	\end{rmk}

Next, 
 we define the random variable:
\[
 \widetilde{\mathcal{R}}_{(4,0)} : \hh \rightarrow \mathbb{R}
\]
such that
\begin{equation*}\label{def:random variable2}
\widetilde{\mathcal{R}}_{(4,0)}(H):=\log{N_{q}\left( \NN\right) } -(4g-4)\log{q}. 
\end{equation*}

 We obtain the following statistical results on $\NN.$

\begin{thm}\label{thm1.6}
	\begin{itemize}
		\item[(1)]
	If $q$ is fixed and $g\rightarrow\infty,$ then 
	\[
	\widetilde{\mathcal{R}}_{(4,0)}
	+\delta_{\gamma/2}\log(1- 1/q^{2})
		\]
	converges weakly to a random variable $\R$  
	such that the characteristic function of $\R$ is given by 
	\begin{eqnarray*}
		\phi_{\R}(t)&=& 
		1
		+
		 \sum\limits_{n=1}\limits^{\infty}
		 \frac{1}{2^{n}n!}
		 \sum\limits_{\substack{P_{1},...,P_{n}     }}
		 \prod\limits_{j=1}\limits^{n}
		 \left(
		 \frac{(1-\mid P_{j} \mid^{-1})^{-it}
		 	+
		 	 (1+\mid P_{j} \mid^{-1})^{-it}-2}{(1+ \mid P_{j}\mid ^{-1})}
	 	 \right),  \, 
	\end{eqnarray*}
 where the inner sum is over distinct monic, irreducible polynomials $P_{1},...,P_{n}$ in $\F_{q^2}[x].$
	
	\item[(2)]
	 If both $q, g\rightarrow\infty,$ then $q\widetilde{\mathcal{R}}_{(4,0)}$ has a standard Gaussian distribution. 
	
\end{itemize}
\end{thm}	
	
	In connection with both in algebraic geometry, integrable system, number theory and in the theory of automorphic forms, there is another interesting moduli space to study and that is the moduli space of Higgs bundles. Whether it is in connection with the space of all solutions of the self-dual equations modulo gauge equivalence on Riemann surfaces \cite{Hitchin}, or giving a geometric interpretation of the theory of elliptic endoscopy which plays crucial role in the proof of the fundamental lemma for unitary groups\cite{Ngo}, \cite{LaNgo}, or to compute the Betti numbers and the Poinca\'re polynomial\cite{HRF08}, \cite{Moz12}, \cite{Sch}, \cite{Mellit20}, the moduli space of stable Higgs bundles always been an active area of research for a decades now. 
	
Assume that $d$ is an odd integer. Let $\text{Higgs}_{2,d}(X)$ stands for the moduli space of stable Higgs bundles of rank 2 and degree $d$ over $X.$ This is a smooth quasi-projective variety defined over $\F_{q}.$ In a similar approach as in Theorem\ref{thm1.3}
	first we give an estimate of the quantity $N_{q}\left( \text{Higgs}_{2,d}(X)\right)$ in terms of $g, q$ and $N.$

	\begin{thm}\label{thm1.1}
		Let $d$ be an odd integer and $X$ be a Galois curve of degree $N$ of genus $g \ge 2$ over $\F_q$. Then there exists a constant $C=C(2,d)$ such that for char$(\F_{q})>C,$ we have
		\begin{align*}
		\log N_{q} \left(\text{Higgs}_{2,d}(X)\right)&=(8g-6)\log q 
		+
		O_{N}\left( 
	 3
	 +q^{-1/2}
	 +\log\left(\frac{\log g}{\log q}\right)  
		 +\frac{1}{q^{2g-2}}\left(\frac{\log g}{\log q}\right)^{N-1}\right).
		\end{align*}
		\end{thm}
	
	Next over the probability space $\H_{\gamma,q}$, we define 
	the random variable 
	\[
	\mathcal{R}^{\a} : \H_{\gamma,q} \rightarrow \mathbb{R}
	\]
	given by
	\[
\mathcal{R}^{\a}(H):= \vw- (8g-6)\log{q}. 
	\]
	Similarly as in Theorem\ref{thm1.4}, we obtain the following statistical results on $\a_{2, d}(H).$
	
	\begin{thm}\label{thm1.2}
		There exists an absolute constant $c>0$ such that 
	the random variable $ \mathcal{R}^{\a}$ has the following decomposition: 
	\[
		 \mathcal{R}^{\a} -C_{q}(2)+\delta_{\gamma/2}\log{\left(1-1/q \right)} = \R_{(\gamma, q)}^{(0)} + \R_{(\gamma, q)}^{(1)}
		+
		O
		\left( q^{-cg}\right), 
		 \]
	 satisfying the following: 
	 \begin{itemize}
	 	\item[(1)] 
	 for a fixed $q,$ and $g \rightarrow\infty,$  $\R_{(\gamma, q)}^{(k)}$ converges weakly to $\R_{k}$
	for $k=0,1,$ where the characteristic function 
	$\phi_{\R_k}(t)$ of $\R_{k}$ is as defined in Theorem\ref{thm1.4}(1).
		
		\item[(2)]
		For a fixed $q,$ we obtain
		\[
		\lim\limits_{\gamma \rightarrow \infty}
		Cov
		\left( 
		\R_{(\gamma, q)}^{(0)}, \R_{(\gamma, q)}^{(1)}
		\right) 
		=\frac{1}{q^{2}} + O\left( \frac{1}{q^{3}}\right).
		\]
		
		\item[(3)]
		If both $q, g\rightarrow\infty,$ then 
		the random variable
		${q^{\frac{(2k+1)}{2}}\R_{k}}$ has a standard Gaussian distribution.
		\end{itemize}	
\end{thm}
	
	We remark here that, all the results on distribution presented in this article  are over the family of hyperelliptic curves $(N=2)$. 
	So, it will be interesting to see whether the proof of Theorem \ref{thm1.4},\ref{thm1.5}, \ref{thm1.6}, and \ref{thm1.2} discussed in section \ref{distribution M(n,d)}, \ref{distribution Ms}, \ref{distribution N tilde}, and \ref{distribution Higgs} respectively can be generalized to a more general set up, that is to study the distribution over a family of non-hyperelliptic curves $(N\geq 3)$.

The layout of the paper is as follows: In the first part of section \ref{preli}, we recall the definitions of zeta function over curves and the Artin $L$-series over function fields. In the second part we recall some basic properties of vector bundles, parabolic vector bundles, Higgs bundles and  moduli spaces. We briefly describe the $\F_{q}$-rational points over these moduli spaces. In Section \ref{distribution M(n,d)}, using induction on the rank of vector bundles, first we give a bound on the number of isomorphism classes of unstable vector bundles (see Proposition \ref{Pdom5.1}), and using this we prove Theorems \ref{thm1.3} and \ref{thm1.4}. The proofs of Theorem \ref{thm1.5} and Theorem \ref{thm1.6} are in Section \ref{distribution Ms} and \ref{distribution N tilde} respectively. In the last Section \ref{distribution Higgs} we discuss Theorem \ref{thm1.1}, and \ref{thm1.2}. We end this section by defining the important notations used in this article.\\  

\noindent
\textbf{Notation:} The notation $f(y)=O(g(y)),$ or equivalently, $f(y)\ll g(y)$ for a non-negative function $g(y)$ implies that there is a constant $c$ such that $|f(y)|\leq cg(y)$ as $y\rightarrow \infty.$ The notation $f(y)=o(g(y))$ is used to denote that $\frac{f(y)}{g(y)} \rightarrow 0$ as $y \rightarrow \infty.$ We use the notation $\mathbb{G}_{m}$ to denote the multiplicative group. 

\section{Preliminaries}\label{preli}
In this section we quickly recall some basic definitions and record some results which will be used later.

\subsection{Zeta functions of curves}   
\label{zeta function basics}
Let $\F_{q}$ be a finite field with $q$ elements and $\overline{\F}_{q}$ be its algebraic closure. 
Let $X$ be a smooth projective geometrically irreducible curve of genus $g \ge 1$ over $\F_{q}$ and $ \bar{X}\,=\,X \times_{\F_{q}} \overline{\F}_{q}$.

Given any integer $r \,>\,0$, let $\F_{q^r}\,\subset \, \overline{\F}_{q}$ be the unique field extension (upto isomorphism) 
of degree $r$ over $\F_{q}$. Let $N_r$ be the cardinality of the set of $\F_{q^r}$- 
rational points 
of $\overline{X}$. Recall that the zeta function of $X$ is defined by 
\begin{equation}
Z_{X}(t)\,=\,\exp\left(\sum_{r >0} \frac{N_rt^r}{r}\right).
\end{equation}
By the Weil conjectures it follows that the zeta function has the form 
\begin{equation}\label{zeta}
Z_{X}(t)\,=\,\dfrac{\prod\limits_{l=1}\limits^{2g}\left(1-\sqrt{q}e(\theta_{l,X})t\right)}{(1-t)(1-qt)},
\end{equation}
where $e(\theta):=e^{2\pi i\theta}.$

Further assume that $X$ is a Galois cover of $\mathbb{P}^1$ with Galois group $G=Gal(R/K)$ of order $N$, where $R:=\F_q(X)$ is the function field of $X$ over $\F_q$ and $K:=\F_q(x)$, the rational function field. For a prime
$\mathfrak{P}$ of $R$, the norm denoted by $|\mathfrak{P}|$ is the cardinality of the residue field of $R$ at $\mathfrak{P}$. The zeta
function $\zeta_R(s)$ is defined by 
\begin{equation}\label{defzeta1.1}
\zeta_R(s)=\prod\limits_{\mathfrak{P} \in R} \left( 1-\frac{1}{|\mathfrak{P}|^s}\right)^{-1}.
\end{equation}
We know that, for the rational function field $K$, prime ideals are in a one-one correspondence with the prime ideals in the polynomial ring $\F_{q}[x],$ only with one exception, and that is the prime at infinity say $P_{\infty}.$ Here $P_{\infty}$ is the discrete valuation ring generated by $\frac{1}{x}$ in $\F_{q}[\frac{1}{x}]$ such that $\deg(P_{\infty})=1.$ By above definition in equation \eqref{defzeta1.1}, the zeta function $\zeta_K(s)$ becomes
\begin{equation*}
\zeta_{K}(s)=\left( 1-\frac{1}{|{P_{\infty}}|^s}\right)^{-1}\prod\limits_{{P}\in \F_{q}[x]} \left( 1-\frac{1}{|{P}|^s}\right)^{-1}=\left( 1-\frac{1}{|{P_{\infty}}|^s}\right)^{-1} \zeta_{\F_{q}[x]}(s).
\end{equation*}
Now $|{P_{\infty}}|=q^{\deg(P_{\infty})}=q.$ And 
\begin{equation*}
\zeta_{\F_{q}[x]}(s)=\sum\limits_{\substack{f \,\,\text{monic}\\ \text{in} \,\, \F_{q}[x]}}\frac{1}{\mid f \mid ^{s}}=\left( 1-q^{1-s}\right)^{-1}. 
\end{equation*}
Therefore,
\[ \zeta_{K}(s)=\left(1-q^{-s}\right)^{-1}\left(1-q^{1-s}\right)^{-1}.
\]
Since $X$ is a smooth projective curve, the zeta function of the curve coincides with
the zeta function of it's function field (see \cite{Ro} for details). More precisely,
\[
Z_X(q^{-s})=\zeta_R(s).
\]
Henceforth we would use $\zeta_R$ and $\zeta_X$ interchangeably to 
denote this zeta function. 
From \eqref{zeta}, we get
\begin{eqnarray}\label{zeta2}
\zeta_{X}(s)=
\frac{\prod\limits_{l=1}\limits^{2g}
	\left(1-\sqrt{q}e(\theta_{l,X})q^{-s}\right)}{(1-q^{-s})(1-q^{1-s})}
=\zeta_K(s)\prod\limits_{l=1}\limits^{2g}\left(1-\sqrt{q}e(\theta_{l,X})q^{-s}\right). 
\end{eqnarray}
Next we recall the Artin L-series for function fields (cf. \cite[Chapter 9]{Ro} for more details). For each prime $P$ of $K$ and a prime $\fP$ of $R$ lying above $P$, 
we denote the Inertia group and the Frobenious element by
$I(\fP/P)$ and $(\fP, R/K)$ respectively.

Let $\rho$ be a representation of $G= \Gal (R/K)$
\begin{center}
	$\rho\,:G\rightarrow\,Aut(V)$
\end{center}
where $V$ is a vector space of dimension $n$ over complex numbers.
Let $\chi$ denotes the character corresponding to $\rho$. For an unramified prime $P$ and $\Re(s)>1$ we define 
the local factor by 
\begin{center}
	$L_P(s,\chi,K)=\det(I -\rho((\fP, R/K))\mid P \mid^{-s})^{-1}$.
\end{center}
Let $\{\alpha_{1}(P), \alpha_{2}(P),....\alpha_{n}(P)\}$ be the eigenvalues of 
$\rho((\fP, R/K)).$ 
In terms of these eigenvalues, we can rewrite the above expression as
\begin{equation}
L_P(s,\chi,K)\,=\prod\limits_{i=1}\limits^{n}
(1-\alpha_{i}(P)\mid P \mid^{-s})^{-1}.
\end{equation}
We note that these eigenvalues $\alpha_{i}(P)$ are all roots of unity 
because $(\fP, R/K)$ is of finite order.

For a ramified prime $P$, the local factor is defined as 
\[
L_P(s,\chi,K)=det(I-\rho((\fP, R/K))_H\mid P\mid^{-s})^{-1}
\]
where $\rho((\fP, R/K))_H$ denote the action of Frobenious automorphism restricted to a subspace $H$ of
$V$ fixed by $I(\fP/P)$.

In either case, we can write 
\begin{equation}\label{local_factor}
L_P(s,\chi,K)=\prod\limits_{i=1}\limits^{n}(1-\alpha_{i}(P)\mid P \mid^{-s})^{-1}.
\end{equation}
where each $\alpha_{i}(P)$ is either roots of unity or zero. 
The Artin $L$-series $L(s,\chi,K)$ or simply $L(s,\chi)$ is defined by
\begin{equation}\label{Artin-factors}
L(s,\chi)= \prod\limits_{P}L_P(s,\chi,K)
\end{equation} 
It is known that if $\rho=\rho_{0}$, 
the trivial representation, then $L(s,\chi)=\zeta_{K}(s),$ and if $\rho=\rho_{reg}$, 
the regular representation, then $L(s,\chi)=\zeta_{R}(s).$

Finally let $\{\chi{_{1}}, \chi_{2}, ... ,\chi_{h}\}$ be the set of irreducible characters 
of the Galois group $G$ with $\chi_{1}= \,\chi_{0},$ the trivial character.
For $i=1,\cdots h$, let $T_{i}=\chi_i(1)$  
be the dimension of the representation space corresponding to 
$\chi_{i}$. 
Then using properties of characters and Artin $L-$series, we get
\begin{equation}\label{l-function-product}
\zeta_{R}(s)\,=\,\zeta_{K}(s)\prod\limits_{i=2}\limits^{h}L(s,\chi_{i})^{T_{i}}  .
\end{equation}

\subsection{Vector bundles over curves} 
\label{vector bundle basics}
Let $X$ be as in subsection \ref{zeta function basics}. 
A vector bundle $E$ on $X$ is a locally free sheaf of $\mathcal O_{X}$-modules of finite rank, where $\mathcal O_{X}$ is the structure sheaf. If $F$ is a subsheaf of a locally free sheaf $E$ for which the quotient $E/F$ is torsion free (and so locally free since $X$ is a curve), then $F$ is called a vector subbundle of $E.$ The rank of the sheaf is denoted by $\text{rank}(E)$. A rank one locally free sheaf is called an invertible sheaf or a line bundle. The degree $\text{deg}(E)$ of a rank $n$ vector bundle $E$ is the degree of it's $n$-th exterior power line bundle 
$\bigwedge^n(E)$ which is also known as determinant line bundle of $E$. For any non-zero vector bundle $E$, the slope of a vector bundle $\mu(E)$ is the rational number $\frac{\text{deg}(E)}{\text{rank}(E)}$. Let $\bar{E}\, =\,E \times_{\F_{q}}\bar{\F}_{q}$ be the extension of $E$ to $\bar{X}$ over $\bar{\F}_{q}$. The vector bundle $\bar{E}$ is called \emph{stable} (resp. \emph{semistable}) if for all proper subbundles $\bar{F}( \neq\, 0, \, \bar{E})$ we have $\mu(\bar{F})\,< \,\mu(\bar{E}) $ (resp. $\mu({\bar{F}})\,\leq \,\mu(\bar{E})$), otherwise it is called \emph{nonsemistable} or \emph{unstable} vector bundle. Over the field $\F_{q},$ a vector bundle $E$ is called stable (resp. semistable) if the corresponding extended vector bundle $\bar{E}$ over $\bar{\F}_{q}$ is stable (resp. semistable). 

For any vector bundle $E$ defined over $\F_{q},$ contains a uniquely determined flag of sub-bundles defined over $\mathbb F_q$ 
\begin{center}
	$0\,=\,F_{0}\,\subsetneqq \,F_{1}\,\subsetneqq...\subsetneqq\, F_{m} \,=\,E$
\end{center}
satisfying following numerical criterion (cf. \cite[Proposition $1.3.9$]{HaNa}). 
\begin{enumerate}
	\item $(F_{i}/F_{i-1})$'s are semistable for $i=1,...,m$.
	\item $\mu(F_{i}/F_{i-1})>\mu(F_{i+1}/F_{i})$ for $i=1,...,m-1.$
\end{enumerate}  
This filtration is often called as Harder-Narasimhan (H-N) filtration or canonical filtration of $E$. The length of the unique flag corresponding to $E$ is called \textit{the length} of $E$ and we denote it by $l(E).$ Also, we use the notation $E^{\vee}$ to denote the dual of the vector bundle $E.$ Assume $X$ is defined over any algebraically closed field (for our purpose $\bar{\mathbb F}_{q}$). For any rational number $\mu,$ let $C(\mu)$ denote the Artinian category of semistable vector bundles on $X$ of slope $\mu.$ For any object $E$ in $C(\mu),$ there is a strictly increasing sequence of vector subbundles  
\begin{center}
	$0\,=\,F_{0}\,\subsetneqq \,F_{1}\,\subsetneqq...\subsetneqq\, F_{m} \,=\,E$
\end{center}
satisfying
\begin{itemize}
	\item $(F_{i}/F_{i-1})$'s are stable for $i=1,...,m$.
	\item $\mu(F_{1})=\mu(F_{2}/F_{1})=...=\mu(F_{m}/F_{m-1})=\mu(E)=\mu.$ 
\end{itemize} 

Such a series is called a \emph{Jordan-H$\ddot{\text{o}}$lder (JH) filtration} of $E.$ The integer $m$ is called the length of the filtration. Though the JH-filtration of a semistable vector bundle $E$ is not unique, the associated grading $gr(E)= \bigoplus \limits_{i} F_{i}/F_{i-1}$ is unique (upto isomorphism). 
Moreover, two semistable bundles $E_{1}$ and $E_{2}$ are called \emph{$S$-equivalent} if $gr E_{1} \cong gr E_{2}$ (cf. \cite{Le}, \cite{Se2}). Note that for stable vector bundles, the $S$-equivalence classes and the isomorphism classes coincides. 

\subsection{Moduli space of vector bundles and it's $\F_{q}$-rational points}
\label{$M(2,0)$ basics}
 Now we consider a line bundle $L$ on $X$ of degree $d$ defined over $\mathbb F_{q}$.   
 Let $M_{L}(r,d)$ (resp. $M^{s}_{L}(r,d)$) be the moduli space of S-equivalence classes of semistable (resp. stable) vector bundles of rank $r$ and determinant isomorphic to $L.$ By going to a finite extension of $\mathbb F_q$ if required, we can assume $M_{L}(r,d)$ is defined over $\mathbb F_q$.
It is well known that $M_{L}(r,d)$ is an irreducible projective variety of dimension $(r^{2}-1)(g-1)$ (cf. \cite{Se2}, \cite{Se4}). 
Further, if $(r,d) \,=\,1,$ then definition of stability and semistability coincides and  $M_{L}(r,d)$ is smooth. 

In this section we will be interested in counting $\F_{q}$-rational points of following three moduli spaces: 
\begin{enumerate}
\item{} Moduli space $M_L(r,d)$ of rank $r$ and degree $d$ vector bundles with fixed determinant $L$, 
when $r$ and $d$ are coprime.
\item{} Moduli space $M^s_{\mathcal O_X}(2,0)$ of rank $2$ and degree $0$ stable vector bundles with fixed determinant $\mathcal O_X$. 
\item{} Sehsadri desingularization $\N$ of $M_{\mathcal O_X}(2,0)$. 
\end{enumerate}

To compute $\F_{q}$-rational points of these moduli spaces, following theorem which is known as {\it Siegel's formula} (see \cite[Section $2.3$]{HaNa}, \cite[Proposition $1.1$]{DeRa}, \cite{GhLe} ) will be used quite frequently:

\begin{thm*}[Siegel's formula]	\begin{equation}\label{Si2.1}
	\sum\limits_{E\in\cM} \dfrac{1}{N_{q}(\Aut E)}\,=\, \frac{1}{q-1}q^{(r^2-1)(g-1)}\zeta_{X}(2)\zeta_{X}(3)...\zeta_{X}(r).
	\end{equation}
where $\mathcal{M}_{L}(r,d)$ denotes the set of all isomorphism classes of rank $r$ vector bundles on $X$ defined over $\mathbb F_q$ with $\det(E) \cong L$, and $\Aut E$ denotes the group scheme of automorphisms of $E$ defined over $\F_{q}$.  
\end{thm*}

\vspace{.5 cm}

\noindent
{\bf (1) $\F_{q}-$Rational points of $M_{L}(r,d)$}, when $gcd(r,d)\,=\,1$: \\
\noindent
By the result of Harder-Narasimhan \cite[Proposition $1.2.1$]{HaNa}, it is known that the number of $\F_{q}$-rational points of $M_{L}(r,d)$ has a bijection with number of isomorphism classes of stable vector bundles over $\F_{q}.$
Also, when $E$ is a stable vector bundle over $\F_{q},$ we know that $\Aut(E) \simeq \mathbb{G}_{m}$ for  (cf. \cite{Ne}), where $\mathbb{G}_{m}$ denote the multiplicative group of $\F_{q}$. Therefore, the Siegel's formula in \eqref{Si2.1} asserts that, 
\begin{equation}\label{Sigen2.2}
N_{q}(\M)\,=\, q^{(r^2-1)(g-1)}\zeta_{X}(2)\zeta_{X}(3)...\zeta_{X}(r)\,-\,\sum\limits_{E\in\Mus} \frac{q-1}{N_{q}(\Aut(E))}.                        
\end{equation} 
where $\Mus\,$ denotes the set of all isomorphism classes of unstable vector bundles. 
For any $E$ in $\Mus,$ it admits a unique Harder-Narasimhan filtration (cf. Section \ref{vector bundle basics}) by subbundles (again defined over $\mathbb F_q$) 
\[ 0\, =\,E_{0}\, \subsetneqq E_{1}\,\subsetneqq...\,\subsetneqq E_{m}\,=\,E.\]
\noindent
We denote the numbers as,
$d_{i}(E):=\deg(E_{i}/E_{i-1}),\, r_{i}(E):=rank(E_{i}/E_{i-1})$ and $\mu_{i}(E):=\mu(E_{i}/E_{i-1}).$
Let
\[HN(n_{1},n_{2},...,n_{m}):=\left\lbrace E \in \Mus\,\vert \,l(E)=m, \, \text{and}\, r_{i}(E)=n_{i}\, \text{for}\, i=1, 2, ...,m\right\rbrace ,\]
\noindent
and 
\begin{equation}\label{unschain}
C_{L}(n_{1},n_{2},...,n_{m})\,:=\, \sum\limits_{E \in HN(n_{1},n_{2},...,n_{m})}\A.
\end{equation}
Then we see that,
\begin{equation}\label{Palpha}
\sum\limits_{E\in\Mus}\frac{1}{N_{q}(\Aut E)}\,=\, \sum\limits_{\substack{(n_{1},n_{2},...,n_{m})\\\sum\limits_{i=1}^{m}n_{i}=r,\, m\geq 2}} C_{L}(n_{1},n_{2},...,n_{m}).
\end{equation}
\noindent
Suppose 
\begin{equation}\label{definition beta}
\beta_{L}(r,d)\,:=\,\sum\A,
\end{equation}
 where the summation extends over isomorphism classes of semistable vector bundles $E$ on $X$ defined over $\mathbb F_q$ of rank $r$ with determinant $L$ of degree $d$. 
We recall the following proposition\cite[Proposition $1.7$]{DeRa}, which will be used later.
\begin{prop}\label{Pbeta}
	\begin{itemize}
		\item[(i)]$\beta_{L}(r,d)$ is independent of $L$ and hence may be written simply as $\beta(r,d).$
		\item[(ii)]\begin{equation}
		C_{L}(n_{1},n_{2},...,n_{m})\,=\, \sum\frac{(N_{q}(J_{X}))^{m-1}}{q^{^\chi \(\tiny{
				\begin{matrix}
				n_{1} & n_{2} &...& n_{m}\\
				d_{1} & d_{2} &...& d_{m}  
				\end{matrix}
			}\)}}\prod\limits_{i=1}^{m}\beta(n_{i},d_{i})
		\end{equation}
		\noindent
		where the summation extends over $(d_{1},d_{2},...,d_{m})\in \mathbb{Z}^{m}$ with $\sum\limits_{i=1}^{m}d_{i}=d$ and $\frac{d_{1}}{n_{1}}>\frac{d_{2}}{n_{2}}>...>\frac{d_{m}}{n_{m}}.$\\
		\noindent
		Here,
		\begin{equation*}
		\chi \(\tiny{
			\begin{matrix}
			n_{1} & n_{2} &...& n_{m}\\
			d_{1} & d_{2} &...& d_{m}  
			\end{matrix}
		}\)=\sum\limits_{i<j}(d_{i}n_{j}-d_{j}n_{i})+\sum\limits_{i<j}n_{i}n_{j}(1-g).
		\end{equation*}
	\end{itemize}
\end{prop}

We complete this subsection by computing the individual terms on the right hand side of \eqref{Palpha} for $r=3.$ 
These computations will be used as the first step in the induction hypothesis to prove Proposition \ref{Pdom5.1}, where we give an asymptotic bound for the number of $\F_{q}$-rational points of isomorphism classes of unstable vector bundles of general rank $r$. 

\begin{prop}\label{unstable rank 3}
	 With the notation as above,
	 \begin{equation}\label{uns3}
	 C_{L}(1,1,1)\leq \frac{q^5\left( N_{q}(J_{X})\right)^{2} q^{3(g-1)}}{(q-1)^{3}(q^2-1)(q^3-1)},
	 \end{equation}
and
\begin{equation}\label{uns4.1}
C_{L}(2,1)=C_{L}(1,2) 
\leq
 \frac{q^{6}N_{q}(J_{X})q^{2(g-1)}}{(q-1)(q^6-1)}\left\lbrace \frac{2q^{3(g-1)}\zeta_{X}(2)}{(q-1)}-\frac{q^{g-1}N_{q}(J_{X})}{(q-1)^{3}(q+1)}-\frac{q^{g}N_{q}(J_{X})}{(q-1)^{3}(q+1)}\right\rbrace .
\end{equation}
\end{prop}

\begin{proof}
Using Proposition \ref{Pbeta}, we see that 
\begin{equation}\label{uns1}
C_{L}(1,1,1)\,=\, \sum\frac{(N_{q}(J_{X}))^{2}}{q^{\chi \(\tiny{
			\begin{matrix}
			1 & 1 & 1\\
			d_{1} & d_{2} &d_{3}
			\end{matrix}
		}\)}}\prod\limits_{i=1}^{3}\beta(1,d_{i}),
\end{equation}
where the summation extends over $(d_{1},d_{2},d_{3})\in \mathbb{Z}^{3}$ with $\sum\limits_{i=1}^{3}d_{i}=d$ and $d_{1}> d_{2} >d_{3}.$ Also,
\begin{equation*}
\chi \(\tiny{
	\begin{matrix}
	1 & 1 & 1\\
	d_{1} & d_{2} & d_{3}  
	\end{matrix}
}\)=2(d_{1}-d_{3}) + 3(1-g).
\end{equation*}
Since $\beta(1, d_{i})=\frac{1}{q-1}$ for $i=1,2,3,$ from \eqref{uns1}, we have 
\begin{equation}\label{uns2}
C_{L}(1,1,1)= 
\frac{
	\left( 
	N_{q}(J_{X})
	\right)^{2} 
	q^{3(g-1)}}{(q-1)^{3}}
\sum
\limits_{      d_{1}>d_{2}>d_{3}        }
\frac{1}{q^{2(d_{1}-d_{3})}
}.
\end{equation}
\noindent
Putting $d_{3}=d-d_{1}-d_{2}$, we see that,
\begin{eqnarray*}
	\sum\limits_{d_{1}>d_{2}>d_{3}}\frac{1}{q^{2(d_{1}-d_{3})}} &=& q^{2d}\sum\limits_{d_{1}>d_{2}>d-d_{1}-d_{2}}\frac{1}{q^{4d_{1}+2d_{2}}}\\
	&=& q^{2d}\sum\limits_{d_{1}>d_{2}}\frac{1}{q^{4d_{1}}}
	\sum\limits_{d_{2}\geq \left[ \frac{d-d_{1}}{2}\right] +1   } \frac{1}{q^{2d_{2}}}.
	\end{eqnarray*}
Note that, in the second summation on the right hand side, either $d_{1}\leq d,$ or $d_{1}>d.$ In either of the two cases, it can be shown that
\begin{equation*}
\sum\limits_{d_{2}\geq \left[ \frac{d-d_{1}}{2}\right] +1   } \frac{1}{q^{2d_{2}}}
\leq \frac{q^2}{q^{d-d_{1}}  (q^2-1)  }.
\end{equation*} 
Therefore, we obtain
\begin{align*}
\sum\limits_{d_{1}>d_{2}>d_{3}}\frac{1}{q^{2(d_{1}-d_{3})}} 
&\leq
\frac{q^{d+2}}{(q^2-1)}
\sum\limits_{d_{1}\geq \left[ \frac{d}{3} \right] +1 }\frac{1}{q^{3d_{1}}}
\leq
\frac{q^5}{(q^2-1)(q^3-1)}.
\end{align*} 
 
Using this in \eqref{uns2}, we get \eqref{uns3}.

Now, let $E$ be in $HN(2,1).$ Therefore, 
we have the Harder-Narasimhan filtration $ 0 \subsetneqq E_{1}  \subsetneqq E $, where $ r(E_{1})=2, r(E/E_{1})=1,$ and the subbundles $E_{1},$ and $E/E_{1}$ are semistable. Assume that $\deg(E_{1})= d_{1},$ and $\det(E_{1})=L_{1}.$ Therefore, $\deg(E/E_{1})=d-d_{1}$ and 
$\det(E/E_{1})=L \otimes L_{1}^{-1}.$ 
The Euler characteristic of $(E/E_{1}^\vee \otimes E_{1}),$ that is $ \chi(E/E_{1}^\vee \otimes E)\,=\, \chi \(\tiny{
	\begin{matrix}
	2 & 1\\
	d_{1} &d-d_{1}
	\end{matrix}
}\)=3d_{1}-2d + 2(1-g)$. 
Now, using Proposition \ref{Pbeta} we get,
\begin{eqnarray*}
	C_{L}(2,1)&=&
	\sum\limits_{\frac{d_{1}}{2}>d-d_{1}}\frac{N_{q}(J_{X})q^{2(g-1)+2d}}{q^{3d_{1}}}\beta(2,d_{1})\beta(1,d-d_{1})\\
	&=& 
	\frac{N_{q}(J_{X})q^{2(g-1)+2d}}{q-1}\sum\limits_{d_{1}>\frac{2d}{3}}\frac{\beta(2,d_{1})}{q^{3d_{1}}}\\
	&=& 
	\frac{N_{q}(J_{X})q^{2(g-1)+2d}}{q-1}
	\left( 
	\beta(2,0)
	\sum\limits_{k \geq \left[ \frac{d}{3}\right] +1   }
	\frac{1}{q^{6k}} + \beta(2,1)
	\sum\limits_{k\geq \left[ \frac{2d-3}{6}  \right] +1 }
	\frac{1}{q^{3(2k+1)}} 
	\right) . \\
	\end{eqnarray*} 
Therefore, 
\begin{equation}\label{uns4}
C_{L}(2,1)
\leq
\frac{N_{q}(J_{X})q^{2(g-1)}}     {(q-1)\left(1-1/q^6 \right)}
\left( 
\beta(2,0)
+ 
\beta(2,1)
\right). 
\end{equation} 

\noindent
Using \eqref{beta(2,0)} in \eqref{calbeta}, which will be discussed in the next section, we obtain
\begin{equation*}
\beta(2,0)= \frac{q^{3(g-1)}}{(q-1)}\zeta_{X}(2)- \frac{ N_{q}(J_{X})q^{g-1}}{(q-1)^{3}(q+1)}.	
\end{equation*}
\noindent
Similarly, we compute $\beta(2,1)$ (cf. Proposition $2.2$ in \cite{ASA} for detailed computation) and putting the values of $\beta(2,0)$ and $\beta(2,1)$ in \eqref{uns4}, we finally obtain \eqref{uns4.1}.

In a similar approach we compute $C_{L}(1,2)$ and one can see that it is equal to the quantity $C_{L}(2,1)$.\\
\end{proof}

\noindent
{\bf (2) $\mathbb F_q$-rational points of $M_{\mathcal{O}_{X}}^{s}(2,0)$}:

\noindent

Let $\mathcal{M}_{\mathcal{O}_{X}}(2,0) \left( \text{resp.}\,\mathcal{M}^{\text ss}_{\mathcal{O}_{X}}(2,0), \mathcal{M}^{\text us}_{\mathcal{O}_{X}}(2,0)\right) $ be the set of all isomorphism classes of rank $2$ vector bundles (resp.  semistable, unstable vector bundles) defined over $\mathbb F_q$ 
on $X$ with trivial determinant $\mathcal{O}_{X}$. 
For brevity we will denote them by $\mathcal M$ (resp. $\mathcal{M}^{\text {ss}}$, $\mathcal{M}^{\text {us}}$).  

We recall from the Definition \eqref{definition beta},
\begin{equation}\label{def:beta(2,0)}
\beta(2,0)\,=\, \sum\limits_{E\in \mathcal{M}^{\text {ss}} }\frac{1}{N_{q}(\Aut(E))}.
\end{equation}
We define
\begin{equation*}\label{def:beta'(2,0)}
\beta'(2,0)\,:=\, \sum\limits_{E\in \mathcal{M}^{\text {us}}}\frac{1}{N_{q}(\Aut(E))}.
\end{equation*}
From the Siegel formula \eqref{Si2.1}, we have
\begin{equation}\label{calbeta}
\beta(2,0)+ \beta'(2,0)= \dfrac{q^{3g-3}}{q-1}\zeta_{X}(2).
\end{equation}

From \cite[equation $(4)$]{BaSe1}, we obtain

\begin{equation}\label{beta(2,0)}		
	\beta'(2,0)=\frac{ N_{q}(J_{X})q^{g-1}}{(q-1)^{3}(q+1)}. 
	\end{equation}

\noindent
Since the stable bundles over $\mathbb F_q$ admit only scaler automorphisms, and the fact that there is a bijection between isomorphism classes of stable vector bundles with trivial determinant (over $\F_{q}$) and the $\F_{q}$-rational points of $M^{\text s}_{\mathcal{O}_{X}}(2,0)$,  from \eqref{def:beta(2,0)}, we can write
\begin{equation}\label{betass}
\beta(2,0)\,=\, \frac{\u}{(q-1)} \,+\, \sum\limits_{E\in \mathcal{M}^{\text{ss}} \setminus \mathcal{M}^{\text{s}}}\frac{1}{N_{q}(\Aut (E))}.
\end{equation}
\noindent
Now we try to compute the second summation on the right hand side of \eqref{betass}. 
Let $K$ be the Kummer variety which is a $2$-fold symmetric product of the Jacobian $J_{\bar X}$. The closed points of $K$ has a set-theoretic bijection 

\[K \leftrightarrow \left\lbrace \text{Isomorphism classes of vector bundles of the form}\,\, \xi\oplus \xi^{-1} \right\rbrace ,\]
where $\xi$ is a degree zero line bundle over $\bar{X}$.  We have a canonical morphism $\phi:\, J_{\bar X} \, \rightarrow \, K$ defined by $ \xi \mapsto \xi \oplus \xi^{-1}$. 
It is known that $K$ has $2^g$ nodal singularities which we denote by $K_0$ and there is a bijection
\begin{equation}\label{defK_{0}}
K_{0} \leftrightarrow \left\lbrace  \xi \oplus \xi: \xi^{2} \cong\mathcal{O}_{\bar{X}},\, \xi\in J_{\bar{X}}
\right\rbrace.
\end{equation}
Without loss of generalities we can assume all $\xi$ which apears in $K_{0}$ are $\F_{q}$-valued points (since they are finitely many and $K$ is projective).  
We also know that $M_{\mathcal O_{\bar X}}(2,0) \setminus M^{s}_{\mathcal O_{\bar X}}(2,0) \cong K$ (cf. \cite{B1}). Moreover, for $E$ being a strictly semistable  vector bundle on $X$ defined over $\mathbb F_q$ with trivial determinant we have an exact sequence over $\bar{\mathbb F_{q}}$
\begin{equation}\label{extn1}
0\longrightarrow \xi\longrightarrow \bar{E}\longrightarrow \xi^{-1}\longrightarrow 0,
\end{equation}
 where $\bar{E}\,=\,E \times_{X} \bar{X}$ and the line bundle $\xi \in J_{\bar X}$ is uniquely determined. 
 
 From the above discussion, we get a surjective set theoretic map 
 \[\theta : \mathcal M^{ss}\setminus\mathcal M^{s} \rightarrow K\]
 which maps a semistable vector bundle $E$ (defined over $\F_{q}$) with $\bar{E}$ as given in \eqref{extn1}, to $gr(\bar E)\,=\,\xi \oplus \xi^{-1}$. Apriori the object $gr(\bar E)$ is defined over $\bar{\mathbb F}_q$ but since $E$ is defined over $\mathbb F_q$, one can show that (see \cite[\S 3]{BaSe1}) both $\xi$ and $\xi^{-1}$ are defined either over $\mathbb F_q$ or $\mathbb F_{q^2}$. 
  Further the image of $\theta$ has the following stratification: 
 \[
 \theta(\mathcal M^{ss}\setminus\mathcal M^{s}) \,\, =\,\,A  \, \sqcup \, B \,\sqcup\, K_{0},
\]
where $A$ and $B$ are defined as follows:
\begin{equation}\label{defA}
A:= \left\lbrace  \xi \oplus \xi^{-1} \in K\setminus K_{0} : \xi \,\text{and} \,\xi^{-1}\,\text{ are both defined over} \,\mathbb{F}_{q} \right\rbrace \,
\text{and}
\end{equation}
\begin{equation}\label{defB}
B:= \left\lbrace  \xi \oplus \xi^{-1}\in K\setminus K_{0}: \xi \,\text{and} \,\xi^{-1}\,\text{  are both defined over $\mathbb {F}_{q^2} $ but not defined over} \,\mathbb{F}_{q} \right\rbrace.
\end{equation}
Clearly, the cardinality

\begin{equation}\label{AB1}
|A|\,=\, \frac{1}{2}(N_q(J_{\bar{X}}) -2^{2g}), 
\end{equation}
and
\begin{equation}\label{AB2}
|B|\,=\,\frac{1}{2}(N_{q^2}(J_{\bar{X}}-J_0) - N_{q}(J_{\bar{X}}-J_0))
\,=\, \frac{1}{2}(N_{q^2}(J_{\bar{X}})-N_{q}(J_{\bar{X}})).
\end{equation}
\noindent
The second equality in \eqref{AB2} is due to the fact that all closed points of $J_{0}$ are $\F_{q}$-valued points.
We set,
\begin{equation}
\label{definition beta1}
\beta_{1}:=\sum\limits_{\substack{E\in \theta^{-1}(A) }}\frac{1}{N_{q}(\Aut(E))}\,+\, \sum\limits_{\substack{E\in \theta^{-1}(B) }}\frac{1}{N_{q}(\Aut(E))},
\end{equation}
and 
\begin{equation}
\label{definition beta2}
\beta_{2}:=\, \sum\limits_{\substack{E\in \theta^{-1}(K_0) }}\frac{1}{N_{q}(\Aut(E))}.
\end{equation}

We have, 
\begin{equation}\label{def2:beta(2,0)}
\beta(2,0)\,=\, \frac{\u}{(q-1)}+ \beta_{1}+\beta_{2}.
\end{equation}
	Therefore, from \eqref{def2:beta(2,0)} and \eqref{calbeta}, we obtain
\begin{equation}\label{Ms1}
\u= q^{3g-3}\zeta_{X}(2)-\left\lbrace \beta'(2,0)+\beta_{1}+\beta_{2}\right\rbrace (q-1).
\end{equation}

Both the term $\beta_1$ and $\beta_2$ are computed in \cite{BaSe1}. While going through the proof given there we noticed in computation of $\beta_1$, they may have missed some term in consideration. Here, we first compute $\beta_{1}.$ We follow the same method as given in [Proposition 3.6, \cite{BaSe1}]. 

  Suppose $E\,\in\,\theta^{-1}(A)$, then there are two possibilities depending on the fact that the extension \eqref{extn1} being split or non-split. In the non-split case $\text{Aut}(E)\,=\, \mathbb G_m$  (over $\mathbb F_q$) and in the split case $\text{Aut}(E)\,=\, \mathbb G_m\,\times\, \mathbb G_m$  (over $\mathbb F_q$) (see [Lemma 3.3, \cite{BaSe1}). Also note that any  extension of $\xi$ with $\xi^{-1}$ 
 is semistable of degree $0$ and two such extension are isomorphic if and only if they are scaler multiple of each other in the extension space $H^1 (X,\xi^{-2})$ which is of dimension $g-1$.  Further if $\xi$ and $\xi'$ are not isomorphic then any two vector bundles $E \in H^1 (X,\xi^{-2})$ and $E' \in  H^1 (X,\xi'^{-2})$ are also non-isomorphic. 
All these facts together tells us that, 
  
  \[
  \sum\limits_{E\,\in\,\theta^{-1}(A)}\frac{1}{N_{q}(\Aut(E))}\,=\, \frac{|A|}{(q-1)^2} \,+\, \frac{2|A| N_q(\mathbb P^{g-2})}{(q-1)}.
  \]

  Now suppose $E \in \theta^{-1}(B)$. Then we know $ E \otimes _{\mathbb F_q} {F_{q^2}} \,=\, \xi \oplus \xi^{-1}$, where $\xi \in J_{\bar{X}} -J_0$ defined over $\mathbb F_{q^2}$. In this case $\text{Aut}(E) \,=\, \mathbb G_m $ over $\mathbb F_{q^2}$ [Lemma 3.5,\cite{BaSe1}]. Clearly if $E_1$ and $E_2$ are in $\theta^{-1}(B)$ and are not isomorphic over $\mathbb F_q,$ then they split non-ismorphically over $\mathbb F_{q^2}$. Hence, 
  
  \[
  \sum\limits_{E\,\in\,\theta^{-1}(B)}\frac{1}{N_{q}(\Aut(E))}\,=\, \frac{|B|}{(q^2-1)}
  \]
  
  Hence we have,
  \begin{equation}\label{Beta1}
  \beta_1\,=\,  \frac{|A|}{(q-1)^2} \,+\, \frac{2|A| N_q(\mathbb P^{g-2})}{(q-1)} \,+\, \frac{|B|}{(q^2-1)}.
  \end{equation}
  
  Now, from Proposition 3.1 in \cite{BaSe1}, we have
 \begin{equation}\label{Beta2}
 	\beta_{2}\,=\,\frac{2^{2g}}{N_{q}(GL(2))}+ \frac{2^{2g}N_{q}(\mathbb{P}^{g-1})}{q(q-1)}.
 	\end{equation}

Putting together the results from \eqref{beta(2,0)}, \eqref{Beta1}, and \eqref{Beta2} along with the size of $A$ and $B$ as in \eqref{AB1} and \eqref{AB2} in \eqref{Ms1} we obtain the following estimate of $\u.$
\begin{prop}\label{Ms2.1}
	The number of $\F_{q}$-rational points of the moduli space of stable bundles $M_{\mathcal{O}_{X}}^{s}(2,0)$ is given by the following expression:
	\begin{align*}\label{Ms}
	\u= q^{3g-3}\zeta_{X}(2)-\frac{\left( q^{g+1}-q^2+q\right) }{(q-1)^2(q+1)}N_{q}(J_{X})-\frac{1}{2(q+1)}N_{q^2}(J_{X})+ \frac{1}{2(q+1)}2^{2g}.
	\end{align*}
\end{prop}

\noindent
{\bf (3) $\mathbb F_q$-rational points of $\N$}:\\

\noindent
 We give a brief description of the Seshadri desingularisation model $\N$ following Seshadri \cite{Se}, and Balaji-Seshadri \cite{BaSe1}.

Let $X$ be any smooth projective curve of genus $\geq 3$ defined over any algebraically closed field $k$ of charecteristic not equal to $2.$ It is known that the moduli space $M_{\mathcal{O}_{X}}(2,0)$ is a projective variety with singular locus being strictly semistable locus which corresponds to  $S$-equivalence classes of semistable vector bundles as explained in previous section. For any vector bundle $V,$ the notion of \emph{parabolic} structure was introduced by Mehta-Seshadri \cite{MeSe}. A parabolic vector bundle can be thought of as a vector bundle with flag structure at each fiber of finitely many points. Let $PV_{4}$ (resp. $PV_{4}^{s}$) denote the category of rank $4$ semistable (resp. stable) parabolic vector bundles $(V, \Delta, \alpha_{\ast})$ with the parabolic structure at a fixed point $p$ in $X$ of flag type $\Delta=(4,3)$ along with parabolic weight $\alpha_{\ast}=(\alpha_{1}, \alpha_{2})$ such that the underlying vector bundle has trivial determinant. 
The weight $\alpha_{\ast}$ can be chosen sufficiently small so that in $PV_{4},$ parabolic semistable $\leftrightarrow$  parabolic stable.
Let $\N$ be the isomorphism classes of parabolic stable vector bundles in $PV_{4}$ such that the endomorphism algebra $\End(V)$ is a specialization of matrix algebra $M_{2\times 2}(k).$ It is known that $\N$ has a structure of a smooth projective variety and there is a morphism 
\[
\pi: \N \rightarrow M_{\mathcal{O}_{X}}(2,0)
\]
which is an isomorphism over $M_{\mathcal{O}_{X}}^{s}(2,0)$ (cf. \cite[Theorem 2.1]{BaSe1}). When $k=\bar{\F}_{q},$ the $\F_{q}$-rational points of $\N$ has been computed in \cite{BaSe1}, which we recall it here.

\begin{prop}(\cite[Theorem $4.2$ ]{BaSe1})
	We have
\begin{equation}\label{pointN2.2}
N_{q}({\N}) \,=\, \u + N_{q}(Y) + 2^{2g}N_{q}(R)+ 2^{2g}N_{q}(S),
\end{equation}
where 
$R$ is a vector bundle of rank $(g-2)$ over $G(2,g),$ the Grassmanian of $2$ dimensional subspaces of $g$ dimensional vector space,
and $S$ is isomorphic to $G(3,g),$ and
$Y$ is a $\mathbb{P}^{g-2}\times \mathbb{P}^{g-2}$ bundle over $K\setminus K_{0},$ 
such that
\begin{equation}\label{pointY}
N_{q}(Y) \,=\, |A| N_{q}(\mathbb{P}^{g-2}\times \mathbb{P}^{g-2}) + |B|N_{q^{2}}(\mathbb{P}^{g-2}),
\end{equation}
 where $A$ and $B$ are the set defined as in \eqref{defA} and \eqref{defB} respectively.
\end{prop}

\subsection{Moduli space of stable Higgs bundles of rank 2} \label{counting point}
Let $X$ be a smooth projective, geometrically connected curve of genus $g,$ defined over the finite field $\F_{q}.$ Let $K_{X}$ be the canonical line bundle of $X.$  A Higgs bundle of rank $r$ and degree $d$ is a pair $(E, \theta)$ with $E$ being a vector bundle of rank $r$ and degree $d$ and $\theta \in \Hom(E, E\otimes K_{X}).$ A Higgs subbundle of $(E, \theta)$ is defined as a subbundle $F\subseteq E$ such that $\theta(F)\subseteq F\otimes K_{X}.$ 
Like in the case of vector bundles, a Higgs bundle $(E, \theta)$ is called \textit{stable} (resp. \textit{semi-stable}) if for all proper Higgs subbundle $F,$ we have $\mu(F)<\mu(E)$ (resp. $\mu(F)\leq \mu(E)$). A Higgs subbundle $F \subset E$ satisfying $\mu(F)>\mu(E)$ is called \textit{destabilizing}. 
When $gcd(r, d)=1,$ the notion of semi-stability and stability both coincide. In particular, for an odd integer $d,$ we consider $\text{Higgs}_{2,d}(X)$ to be the moduli space of stable Higgs bundles over $X$, which is a smooth, quasi-projective, cohomologically pure variety of dimension $(8g-6).$ Without loss of generality we can assume both $X$ and $\text{Higgs}_{2,d}(X)$ are defined over $\F_{q}.$ 
We denote the cardinality of the set of $\F_{q}$-rational points over $\text{Higgs}_{2,d}^{st}(X)$ by $N_{q}(\b).$ 
The computation of the quantity $N_{q}(\b)$ is a two fold analysis due to Schiffmann \cite{Sch}. First step is to relate the stable Higgs bundles in terms of geometrically indecomposible vector bundles on $X,$ and then expressing the number of such geometrically indecomposible vector bundles in the form of a unique rational function in Weil numbers associated to the curve $X$.

We state \cite[Theorem 1.2]{Sch}, which is more general and applicable for any rank $r$ and degree $d$ with the condition $gcd(r,d)=1.$ For our purpose, here we state the theorem only for rank $2$ cases.
\begin{thm}[Schiffmann]\label{Schiffman 1.1}
	Let $d$ be an odd integer. 
	There exists an explicit constant $C=C(2,d)$ such that for any smooth projective geometrically connected curve $X$ of genus $g$ defined over $\F_{q}$ with char$(\F_{q})>C$, we have
	\[
	N_{q}(\b)=q^{(4g-3)}A_{g,2,d}
	\]
	where $A_{g,2,d}$ denote the number of geometrically indecomposable vector bundles on $X$ of rank $2$ and degree $d.$ 
\end{thm}
Next we state \cite[Theorem 1.6]{Sch} for rank 2.
\begin{thm}[Schiffmann]\label{Schiffman 1.6}
	Suppose the eigenvalues of the Frobenious acting on the first cohomology of  any smooth projective geometrically connected curve $X$ of genus $g$ defined over $\F_{q}$ are $\alpha_{1},...,\alpha_{2g}$ with $\alpha_{i+g}=q\alpha_{i}^{-1}$ for $i=1,...,g.$ Then for any odd integer $d,$ the number of geometrically indecomposable vector bundles on $X$ of rank $2$ and degree $d,$ that is
	\begin{equation}\label{def:Ag2}
	A_{g,2,d}\,=\,A_1+A_{2}+A_{3}
	\end{equation}
	where
	\begin{align*}
	A_{1}&:=\frac{   \prod\limits_{i=1}^{2g}(1-\alpha_{i})   \prod\limits_{i=1}^{2g}(1-q\alpha_{i})       }  
	{     (q-1)(q^2-1)    }\\
	A_{2}&:=-
	\frac{   \prod\limits_{i=1}^{2g}(1-\alpha_{i})   \prod\limits_{i=1}^{2g}(1+\alpha_{i})       }  
	{     4(q+1)   }\\
	A_{3}&:=
	\frac{   \prod\limits_{i=1}^{2g}(1-\alpha_{i})^{2}       }  
	{     2(q-1)    }\left( 
	\frac{1}{2} -\frac{1}{q-1}
	-
	\sum\limits_{i=1}^{2g}
	\frac{1} {(1-\alpha_{i})    }
	\right). 
	\end{align*}
\end{thm}
Note that, the quantity $A_{g,2,d}$ is come out to be independent of the degree of the vector bundles and we can drop the dependence on $d$ from the notation and we simply denote it by $A_{g, 2}$ from now on.
In fact, for general rank and degree, the quantity $A_{g, r, d}$ does not depend on the degree $d$. Although, it is not very clear from \cite[Theorem 1.6 ]{Sch}, it is evident from \cite[Theorem 1.1]{Mellit20}.

\section{Distribution on $M_{L}(r,d)$ } \label{distribution M(n,d)}

Let $X$ be a Galois curve of degree $N$ and genus $g \geq 2$ over a finite field $\F_{q}$ and $M_{L}(r,d)$ be the moduli space of stable vector bundles on $X$ of rank $r$ with fixed determinant $L$ of degree $d,$ such that $gcd(r,d)=1.$ 

\subsection{Proof of Theorem \ref{thm1.3}}
With the same set up as in subsection \ref{zeta function basics}, we recall that, for the smooth projective Galois curve $X,$ $R:=\F_q(X)$ is a geometric Galois extension 
of $\F_q(x)$ with Galois
group $G$ of order $N$.  From equation (\ref{zeta2}), and (\ref{l-function-product}), we get
\[
\prod\limits_{i=2}\limits^{h}L(s,\chi_{i})^{T_{i}}
=\prod\limits_{l=1}\limits^{2g}(1-\sqrt{q}e(\theta_{l,X})q^{-s}).
\]
Using (\ref{local_factor}) and (\ref{Artin-factors}) on the left hand side, we get
\begin{equation}\label{local-Artin5.1}
\prod_{i=2}^h \prod\limits_{P}\prod\limits_{j=1}\limits^{T_{i}}(1-\alpha_{i,j}(P)\mid P\mid^{-s})^{-T_i}
=\prod\limits_{l=1}\limits^{2g}(1-\sqrt{q}e(\theta_{l,X})q^{-s})
\end{equation}
where the product on the left hand side is over all monic, irreducible polynomials $P\;\text{in}\; \F_{q}[x],$
and $P_\infty$. 
These $\alpha_{i,j}(P)$'s are either roots of unity or zero.

 Now, Taking logarithms on both sides of \eqref{local-Artin5.1} and equating the coefficient of $q^{-ms}$ 
for any positive integer $m$, we obtain 
\begin{equation}\label{theta-lambda-equality5.1}
-q^{m/2}\sum_{l=1}^{2g}e(m\theta_{l,X})=\sum\limits_{\deg\, f = m}
\Lambda(f)\sum\limits_{i=2}\limits^{h}T_{i}\sum\limits_{j=1}\limits^{T_{i}}\alpha_{i,j}(f)
\end{equation}
where the sum on the right is over all monic irreducible polynomials of degree $m$ over $\F_q$ and $\Lambda $ is the analogue of Von Mangoldt function defined as
\begin{equation*}
\Lambda(f):= \begin{cases}
\deg P &\text{ if } f=P^{k}\text{ for some monic, irreducible } P \in \F_{q}[x]\\  
0  &\text{  otherwise.}
\end{cases}
\end{equation*}
From the definition of zeta function as in \eqref{zeta2}, for any integer $k\geq2$ we can write 
\begin{eqnarray*}
	\zeta_{X}(k)&=& \frac{q^{(2k-1)}\prod\limits_{l=1}\limits^{2g}\(1-q^{-(2k-1)/2}e(\theta_{l,X})\)}{(q^{k}-1)(q^{k-1}-1)}.
\end{eqnarray*}
Taking logarithm on both sides of the above equation we get
\begin{equation}\label{zetadef5.2}
\log \zeta_{X}(k) - (2k-1)\log {q} + \log \left( (q^{k}-1)(q^{k-1}-1)\right)  \,=\, \sum\limits_{l=1}\limits^{2g}\log (1-q^{-(2k-1)/2}e(\theta_{l,X})). 
\end{equation}
For any positive integer $Z$, we define,
\begin{equation}\label{def: epsilon1}
 \epsilon_{1,Z} := -\sum\limits_{m\leq Z}q^{-(2k-1)m/2}m^{-1}\sum\limits_{l=1}\limits^{2g}e(m\theta_{l,X}),
 \end{equation}
and 
\begin{equation}\label{def: epsilon2}
\epsilon_{2,Z} := -\sum\limits_{m\geq Z+1}q^{-(2k-1)m/2}m^{-1}\sum\limits_{l=1}\limits^{2g}e(m\theta_{l,X}).
\end{equation}
Putting these in \eqref{zetadef5.2}, we can write 
\begin{equation}\label{log-zeta5.1}
\log{\zeta_{X}(k)} - (2k-1)\log{q} + \log\left( (q^k-1)(q^{k-1}-1)\right)  \,=\,\epsilon_{1,Z} + \epsilon_{2,Z}.
\end{equation}  
In the next result we estimate $\epsilon_{1,Z}$ and $\epsilon_{2,Z}$.
\begin{lem}\label{epsilon-Z5.1}
	For  $Z \geq 2$, we have
	\[| \epsilon_{1,Z}| \, \leq \,(N-1)\left( \frac{1}{q-1} + \frac{1}{q^{k}} (1.5 + \log{Z}-\log{2})\right).\] 
	\noindent
	and 
	\[| \epsilon_{2,Z}| \, \leq \, \frac{2g}{(Z + 1)} \, \frac{1}{q^{(2k-1)(Z + 1)/2}} \, \frac{1}{(1-q^{-(2k-1)/2})}.\]
	Moreover,
	\[| \epsilon_{1,1}| \, \leq \, (N-1)\left( \frac{1}{q^{k-1}} + \frac{1}{q^{k}}\right) \]
	and
	\[| \epsilon_{2,1}|\, \leq\, \frac{2g}{q^{(2k-1)}-q^{(2k-1)/2}}.\]
\end{lem}

\begin{proof}
	
	Let $Z\,\geq\, 2$. Using (\ref{theta-lambda-equality5.1}) in the definition of $\epsilon_{1, Z}$ in \eqref{def: epsilon1}, we get
	\begin{eqnarray*}
		| \epsilon_{1, Z}|
		&\leq & \sum\limits_{m\leq Z}q^{-k m}m^{-1}\sum\limits_{deg f = m}\Lambda(f)\sum\limits_{i=2}\limits^{h}T_{i}^{2}. 
	\end{eqnarray*}
	Now, using the property that $\sum\limits_{i=2}\limits^{h}T_{i}^{2}=N-1,$ and 
	$\sum\limits_{\deg\, f = m}\Lambda(f)=q^{m}+1,$ we have
	\begin{eqnarray*}
		| \epsilon_{1,Z}| &\leq &  \sum\limits_{m\leq Z}q^{-k m}m^{-1}(q^{m}+1)(N-1) \\
		&=& (N-1)\left( \sum\limits_{m\leq Z}\dfrac{1}{mq^{(k-1)m}} + \sum\limits_{m\leq Z}\dfrac{1}{mq^{km}}\right). 
	\end{eqnarray*}
Also, we know for a positive integer $n,$
\begin{equation*}
\sum\limits_{m \geq 1}\frac{1}{mn^{m}} \leq -\log(1-n^{-1})\leq \frac{1}{n-1}.
\end{equation*}
That implies
\begin{eqnarray*}
	|\epsilon_{1,Z}| &\leq & (N-1)\left( \frac{1}{q^{(k-1)}-1} + \frac{1}{q^{k}}\sum\limits_{m\leq Z}\dfrac{1}{m}\right)  \\
	& \leq & (N-1)\left( \frac{1}{q^{k-1}-1} + \frac{1}{q^{k}}(1.5 + \log{Z} -\log{2})\right).
\end{eqnarray*}
Similarly, from the definition of $\epsilon_{2, Z}$ in \eqref{def: epsilon2}, we have
\begin{eqnarray*}
	| \epsilon_{2,Z}| &= & 
	\bigg| 
	\sum\limits_{m\geq Z + 1}q^{-(2k-1)m/2} m^{-1}\sum\limits_{l=1}\limits^{2g}-e(m\theta_{l,X})
	\bigg| \\
	&\leq & 2g\sum\limits_{m\geq Z+1}q^{-(2k-1)m/2}m^{-1} \\
	&\leq & \frac{2g}{(Z + 1)} \, \frac{1}{q^{(2k-1)(Z + 1)/2}} \, \frac{1}{(1-q^{-(2k-1)/2})}.
\end{eqnarray*}
Also, for $Z=1,$ trivially we get the desired bound for $\epsilon_{1,1},$ and $\epsilon_{2,1},$ and that completes the proof of Lemma \eqref{epsilon-Z5.1}.
\end{proof}

\begin{prop}\label{log-zeta5.2}
	With all the notations as above, for suitable absolute constants $c_{1}>0$ and $c_{2}>0$ and assuming $\log(g)>\kappa \log(Nq)$ for a sufficiently large absolute constant $\kappa >0$ (independent of $N$), we obtain
	\begin{equation*}
	\left|  \log{\zeta_{X}(k)}\right|  \leq \frac{c_{1}N}{\sqrt{q}}+\frac{c_{2}N\log{\log{g}}}{q^{k}}
	\end{equation*}
	for any $k\geq 2.$
\end{prop}
\begin{proof}
	If  $ \frac{2g}{(q^{2k-1}-q^{(2k-1)/2})}\, \geq \, q^{-1/2}(N-1)$, 
	using Lemma \ref{epsilon-Z5.1} with  
	\[Z =\frac{2}{3}\,\frac{\log{\left( \frac{2g\sqrt{q}}{(1-q^{-(2k-1)/2})(N-1)} \right) }}
	{\log{q}}\, \geq 2
	\]
	we get
	\begin{equation*}
	| \epsilon_{1,Z}| \leq (N-1)\left( \frac{1}{q-1} + \frac{1}{q^{k}}(1.5 -\log{2})\right)
	+ \frac{(N-1)}{q^{k}}\log{\left(\frac{2}{3}\,\frac{\log{\left( \frac{2g\sqrt{q}}{(1-q^{-(2k-1)/2})(N-1)} \right) }}
		{\log{q}}\right)}.
	\end{equation*}
	and
	\[
	| \epsilon_{2,Z}|  \leq \frac{(N-1)}{q^{k}}\log\left(\frac{2}{3}\frac{\log{\left( \frac{2g\sqrt{q}}{(1-q^{-(2k-1)/2})(N-1)} \right) }}
	{\log{q}} \right) 
	\]
Note that $\frac{1}{q-1}\leq \frac{1}{q}+\frac{2}{q^2}.$
Hence in this case
\begin{equation}\label{Z2}
| \epsilon_{1,Z}| + | \epsilon_{2,Z}|
\leq \,(N-1)
\Bigg(
\frac{1}{q} +\frac{2}{q^2}+ 
\frac{2}{q^{k}}
\log{
	\Bigg(
	\frac{2\log{
			\left( \frac{2g\sqrt{q}}{(N-1)(1-q^{-(2k-1)/2})}
			 \right) }}       {3\log{q}}
	\Bigg)
} 
+ \frac{1}{q^{k}}
\Bigg). 
\end{equation}
Now, suppose $ \frac{2g}{(q^{2k-1}-q^{(2k-1)/2})}\, < \, q^{-1/2}(N-1)$. 
Again using Lemma \ref{epsilon-Z5.1} with $Z=1,$ we get
\begin{equation}\label{Z1}
| \epsilon_{1,1}|+ | \epsilon_{2,1}| \leq 
(N-1)\left\lbrace \frac{1}{\sqrt{q}}+\frac{1}{q^{k-1}}+\frac{1}{q^{k}}\right\rbrace .
\end{equation}
Using \eqref{Z2}, and \eqref{Z1} in \eqref{log-zeta5.1}, we get 
\begin{align}\label{log-zeta2}
&\mid \log{\zeta_{X}(k)} -(2k-1)\log{q} + \log\{(q^k-1)(q^{k-1}-1)\}\mid \nonumber \\
&\le(N-1)\left( \frac{c}{\sqrt{q}} + 
\frac{2}{q^{k}}\left\lbrace 
\log{\left(\frac{2\log{\left( \frac{2g\sqrt{q}}{(N-1)(1-q^{-(2k-1)/2})} \right) }}{3\log{q}}\right)} \right\rbrace\right)
\end{align}
for some absolute constant $c> 0.$
After simplifying we can write
\begin{eqnarray*}
	\frac{2\log{\left( \frac{2g\sqrt{q}}{(N-1)(1-q^{-(2k-1)/2})} \right) }}{3\log{q}}
	&=& \frac{2\log{g}}{3\log{q}} \left(1 + O\left(\frac{\log{(Nq)}}{\log{g}} \right)  \right) .
\end{eqnarray*}

Therefore,
\begin{align*}
\frac{(N-1)}{q^{k}}\log \left( \frac{2\log{g}}{3\log{q}} \left(1 + O\left(\frac{\log{(Nq)}}{\log{g}} \right)  \right) \right)
&= \frac{(N-1)\log{\log{g}}}{q^{k}} + O\left( \frac{N}{q^{k-1}}
+\frac{N}{q^{k}}\frac{\log{(Nq)}}{\log{g}}\right) 
\end{align*}
\noindent
under the assumption $ \log{g}> \kappa \log{(Nq)}$ for a sufficiently large constant $\kappa >0.$
Also, we see that, 
\begin{equation*}
\log(q^{k}-1)(q^{k-1}-1)
= (2k-1)\log{q}\,+\,O\left( \frac{1}{q^{k-1}}\right) .
\end{equation*}
Using above results in \eqref{log-zeta2}, we obtain the desired bound
%
\[ \left|  \log{\zeta_{X}(k)}\right|  \leq \frac{c_{1}N}{\sqrt{q}}+\frac{c_{2}N\log{\log{g}}}{q^{k}}\]
\noindent
with suitable choice of two constants $c_{1}>0$ and $c_{2}>0$ and hence the Proposition.

\end{proof}

As a direct consequence of Proposition \ref{log-zeta5.2}, we obtain the following result. 
\begin{prop}\label{log-zeta5.3}
	There exists an absolute constant $c'>0$ such that
	\[
	(\log\,g)^{\frac{-c'N}{q^k}}\exp\left(\frac{-c'N}{\sqrt{q}} \right) \leq\zeta_{X}(k)\leq (\log\,g)^{\frac{c'N}{q^k}}\exp\left( \frac{c'N}{\sqrt{q}}\right), 
	\] 
	for any $k\geq 2,$ whenever $\log{g}>\kappa \log{(Nq)}$ for a sufficiently large constant $\kappa>0.$ 
\end{prop}


Now we recall the following result of Xiong and Zaharescu (Theorem 1 of \cite{XiZa}).
\begin{lem}\label{Xiong-Zaha}
	Let $X$ be a Galois curve of degree $N$ of genus $g \geq 1$ over $\F_{q}.$
	Then
	\begin{equation}\label{log-Jacobian}
	\left. \mid\log{\left( N_{q}(J_{X})\right)}-g\log{q}\right. \mid\, 
	\leq\, (N-1)\left(\log{\max\left\lbrace 1,\dfrac{\log{\left(\frac{7g}{(N-1)}\right)}}{\log{q}}\right\rbrace }
	+3\right).
	\end{equation}
\end{lem}
A direct consequence of Lemma\ref{Xiong-Zaha} is as follows.
\begin{prop}\label{log_Jac5.1}
	For $\log g>\kappa_{1} \log q,$ where $\kappa_{1}>0$ is a large constant, we have
	\[
	 q^{g}
	 \left( 
	 \log\frac{7g}{N-1}
	 \right) 
	 ^{-3(N-1)}
	 \leq N_{q}
	 (J_{X})
	 \leq 
	 q^{g}
	 \left( 
	 \log\frac{7g}{N-1}
	 \right) 
	 ^{3(N-1)}.
	 \] 
\end{prop}

Next we recall the definition of $C_{L}(n_{1},n_{2},...,n_{k})$ in equation \eqref{unschain}. 
Based on the above mentioned two results on the bounds of the quantity $N_{q}(J_{X})$ and $\zeta_{X}(k)$,  we get the following result.
\begin{prop}\label{Pdom5.1}
	Let $E$ be an unstable vector bundle in $\Mus.$ Then for any partition $(n_{1},n_{2},...,n_{k})$ of the rank $r\geq 3,$ we have
	\[
	C_{L}(n_{1},n_{2},...,n_{k})
	=o\left(
	q^{(r^2-\frac{r}{c})(g-1)}\zeta_{X}(2)\zeta_{X}(3)...\zeta_{X}(r-1) 
	\right)
	\] 
	for all $k\geq 2,$ and $c>1.$
\end{prop}
\begin{proof}
	We use induction on the rank $r\geq 3.$ 
For $r=3,$
	using Proposition\ref{log-zeta5.3} and Proposition \ref{log_Jac5.1} in \eqref{uns4.1}, we have
	\[
	\lim_ {g \rightarrow \infty}\dfrac{C_{L}(2,1)}{q^{(3^2-\frac{3}{c})(g-1)}\zeta_{X}(2)}=0.
	\]
	Similarly one can show that $C_{L}(1,2)$ and $C_{L}(1,1,1)$ are of the size $o(q^{(3^2-\frac{3}{c})(g-1)}\zeta_{X}(2)).$ Hence the induction hypothesis holds for the initial case.

	Now, suppose the statement is true for all partitions of $m,$ where $m$ is the rank of any vector bundle $E$ in $\mathcal{M}_{L}^{\text{us}}(m,d),$ and $m < r.$ That is for any $m<r,$ 
	\[C_{L}(m_{1},m_{2},...,m_{k})=o\left(q^{(m^2-\frac{m}{c})(g-1)}\zeta_{X}(2)\zeta_{X}(3)...\zeta_{X}(m-1) \right)\]
	\noindent
	for all $k\geq 2$ such that $\sum\limits_{i=1}^{k}m_{i}=m$. We want to show that the statement is true for rank $r.$
Suppose $E$ is in $HN(n_{1},n_{2},...,n_{m}),$ such that $E$ has the H-N filtration 
\[
 0=E_{0} \subsetneqq E_{1}...\subsetneqq E_{m}=E.
 \]
We write $M=E/E_{1}.$ Using Proposition $4.4$ in \cite{NaSe} it can be shown that, there is no nonzero homomorphism from $E_{1}$ to $M.$ Hence every automorphism of $E$ keeps $E_{1}$ invariant and it goes down to an automorphism of the quotient $M.$ Thus we get a well defined map :
\[
\Phi: \Aut E \rightarrow \Aut E_{1} \times \Aut M,
\]
such that, for any $f$ in $\Aut E,$ $\Phi(f)=(f\vert _{E_{1}}, p\circ f)$ where $\xymatrix{E\ar[r]^{f} & E\ar[r]^{p} & M}.$ Note that for any $g \in \Hom(M,E_{1}),$ we have the image $\Phi(Id+g)=(Id\vert_{E_{1}}, Id\vert_{M})$ for $\left( Id + g\right) $ in $\Aut E.$ Therefore, $ Id + H^{0}(\bar{X}, \Hom(M, E_{1}))$ is contained in $\ker(\Phi).$ Conversely let any $f \neq Id,$ is in $\ker(\Phi),$ that is $\Phi(f)=(Id\vert_{E_{1}}, Id\vert_{M}).$ We see that $(f-Id)\vert _{E_{1}}=0$ and $p\circ (f-Id)=0,$ and since $\Hom(E_{1}, M)=0,$ therefore, $(f-Id)$ is in $\Hom(M, E_{1}).$ So any $f$ in $\ker(\Phi),$ is in $Id +H^{0}(\bar{X}, \Hom(M, E_{1})), $ and hence $\ker(\Phi) \cong H^{0}(\bar{X}, \Hom(M, E_{1}))=T$(say).

Next we consider the action of the group $G:=\Aut E_{1} \times \Aut M$ on $ H^{1}(\bar{X}, \Hom(M, E_{1})),$ that is, on equivalence classes of extensions of $M$ by $E_{1}.$ For simplicity we denote $H^{1}(\bar{X}, \Hom(M, E_{1}))$ by $S.$ We denote any equivalene class of short exact extensions 
\[
\xymatrix{0\ar[r] & E_{1}\ar[r]^{\alpha} & E\ar[r]^{\beta} & M\ar[r] & 0}
\]
 in $S$ by $[E; (\alpha, \beta)].$ For any $(\phi_{1}, \phi_{2})$ in $G,$ define the action by 
\[
(\phi_{1}, \phi_{2})\cdot [E; (\alpha, \beta)]= [E; (\alpha\circ \phi_{1}, \phi_{2}\circ\beta)].
\]
We know that, any two extensions $(E; (\alpha_{1}, \beta_{1}))$ and $(E; (\alpha_{2}, \beta_{2}))$ are isomorphic if the following diagram commute:
\begin{center}

	\begin{tikzpicture}[descr/.style={fill=white,inner sep=2.5pt}]
\matrix (m) [matrix of math nodes, row sep=3em,
column sep=3em]
{ 0 & E_{1} & E & M & 0\\
	0& E_{1} & E & M & 0\\ };
\path[->,font=\scriptsize]
(m-1-1) edge node[auto] {} (m-1-2)
(m-1-2) edge node[auto] {$ \alpha_{1} $} (m-1-3)
(m-1-3) edge node[auto] {$ \beta_{1} $} (m-1-4)
(m-1-4) edge node[auto] {} (m-1-5)
(m-1-2)edge node[left]  {$ \phi_{1} $} (m-2-2)
(m-1-3)edge node[left] {$\phi $} (m-2-3)
(m-1-4)edge node[left] {$ \phi_{2} $} (m-2-4)
(m-2-1) edge node[auto] {} (m-2-2)
(m-2-2) edge node[auto] {$ \alpha_{2} $} (m-2-3)
(m-2-3) edge node[auto] {$ \beta_{2} $} (m-2-4)
(m-2-4) edge node[auto] {} (m-2-5);
\end{tikzpicture}
\end{center}
where $\phi_{1}, \phi$ and $\phi_{2}$ are isomorphisms. When $\phi_{1}=Id$ and $\phi_{2}=Id$, we get the equivalence class. From the definition, it is easy to see that, two such extensions in $S$ are isomorphic if and only if they are in the same orbit under this action and the isotropy subgroup of $[E; (\alpha, \beta)]$
denoted by $G_{[E; (\alpha, \beta)]}$is same as the image of $\Aut E$ under the map $\Phi.$ We denote the orbit of $[E; (\alpha, \beta)]$ in $S$ by $G\cdot [E; (\alpha, \beta)].$ Therefore, we have
\begin{equation*}
C_{L}(n_{1},n_{2},...,n_{m})= \sum\limits_{E_{1},M}\,\sum\limits_{[E; (\alpha, \beta)] \in S}\frac{1}{\vert \Aut E\vert \, \vert G \cdot [E; (\alpha, \beta)]\vert}.
\end{equation*}
Also, $\vert G \cdot [E; (\alpha, \beta)]\vert = \left[ G:G_{[E; (\alpha, \beta)]}\right] = \frac{\vert \Aut E_{1} \times \Aut M \vert }{\vert Img(\Phi)\vert}= \frac{\vert \Aut E_{1} \times \Aut M \vert \, \vert \ker(\Phi)\vert}{\vert \Aut E\vert}.$ 
So, 
\begin{equation*}
C_{L}(n_{1},n_{2},...,n_{m})= \sum\limits_{E_{1},M}\frac{\vert S \vert }{\vert \Aut E_{1}\vert \, \vert \Aut M\vert \, \vert T \vert} = \sum\limits_{E_{1},M}\frac{1}{\vert \Aut E_{1}\vert \, \vert \Aut M\vert q^{\chi(M^{\vee}\otimes E_{1})}},
\end{equation*}
where $\chi(M^{\vee}\otimes E_{1})$ denotes the Euler characteristic of $(M^{\vee}\otimes E_{1}).$ The summation extends over all pairs of bundles $(E_{1},M)$ where $E_{1}$ is semistable of rank $n_{1},$ and $M$ has H-N filtration of length $(m-1)$ and has determinant equal to $L \otimes (\det E_{1})^{-1}$ such that $\mu(E_{1})> \mu_{1}(M)>...>\mu_{m-1}(M)$ with $r_{i}(M)=n_{i+1},\, \, i=1, 2,...,m-1.$ 

Now, let $J_{X}^{d_{1}}$ be the variety of isomorphism classes of line bundles of degree $d_{1}$ on $X.$ 
Therefore,
\begin{equation}\label{induction1}
C_{L}(n_{1}, n_{2},...,n_{m})= \sum\limits_{\substack{d_{1} \in \Z\\
		 \frac{d_{1}}{n_{1}}> \frac{d-d_{1}}{r-n_{1}}}} 
	 \sum\limits_{L_{1} \in J_{X}^{d_{1}}(\F_{q})} 
	 \frac{C_{L \otimes L_{1}^{-1}}(n_{2}, n_{3},...,n_{m})}{q^{\chi \(\tiny{
			\begin{matrix}
			n_{1} &n- n_{1}\\
			d_{1} & d-d_{1}
			\end{matrix}
		}\)}} \sum\frac{1}{\mid \Aut E_{1} \mid},
\end{equation}
where the last sum on the right hand side of the above equation is over isomorphism classes of semistable bundles $E_{1}$ of rank $n_{1}$ with $\det(E_{1}) \cong L_{1},$
 which is nothing but $\beta_{L_{1}}(n_{1}, d_{1})$ as defined in \eqref{definition beta}.
Also,
 \[
 \chi\left(M^{\vee}\otimes E_{1} \right) =  \chi \(\tiny{
	\begin{matrix}
	n_{1} &r- n_{1}\\
	d_{1} & d-d_{1}
	\end{matrix}
}\)
=rd_{1}-n_{1}d-n_{1}(r-n_{1})(g-1).
\]
Note that $N_{q}({J_{X}})=N_{q}(J_{X}^{d})$ for any degree $d$ (cf. \cite{DeRa}).
Therefore, for any partition $(n_{1},n_{2},...,n_{k})$ of $r,$ using Proposition \ref{Pbeta}$(i)$ in \eqref{induction1}, we write 
\begin{align}\label{induction2}
	C_{L} (n_{1},n_{2},...,n_{k})
	&= \sum\limits_{\substack{d_{1}\in \Z \\ \frac{d_{1}}{n_{1}}>\frac{d}{r}}}\, \frac{N_{q}(J_{X})q^{n_{1}(r-n_{1})(g-1)+n_{1}d}}{q^{rd_{1}}}\beta(n_{1},d_{1})C_{L\otimes L_{1}^{-1}} (n_{2},...,n_{k}).
\end{align}

Using \eqref{Palpha} in the Siegel's formula \eqref{Si2.1} for rank $n_{1},$ we obtain
\begin{align*}
\beta(n_{1},d_{1})
&=
 \frac{1}{q-1}q^{(n_{1}^2-1)(g-1)}\zeta_{X}(2)\zeta_{X}(3)...\zeta_{X}(n_{1})\\
&-
 \sum\limits_{\substack{(p_{1},p_{2},...,p_{l})\\\sum\limits_{i=1}^{l}p_{i}=n_{1},\, l\geq 2}} C_{L_{1}}(p_{1},p_{2},...,p_{l}).
\end{align*}
Using induction hypothesis for rank $n_{1}<r,$ we have
\begin{align}\nonumber
\beta(n_{1},d_{1})&= \frac{1}{q-1}q^{(n_{1}^2-1)(g-1)}\zeta_{X}(2)\zeta_{X}(3)...\zeta_{X}(n_{1})\\\label{beta-n5.1}
&+
 o\left(S(n_{1})q^{(n_{1}^2-\frac{n_{1}}{c})(g-1)}\zeta_{X}(2)\zeta_{X}(3)...\zeta_{X}(n_{1}-1) \right),
\end{align}
where $S(n_{1})$ is the number of partitions of $n_{1}.$
Also, 
\begin{equation*}
\sum\limits_{\substack{d_{1}\in \Z \\ \frac{d_{1}}{n_{1}}>\frac{d}{r}}}\,\frac{1}{q^{rd_{1}}}\,=\,
O\bigg( 
\frac{1}{q^{dn_{1}}}
\bigg). 
\end{equation*}

Using induction hypothesis on $C_{L\otimes L_{1}^{-1}} (n_{2},...,n_{k})$ for rank $(r-n_{1}),$ and putting the bound from\eqref{beta-n5.1} in \eqref{induction2}, we obtain
\begin{align*}
	C_{L} (n_{1},n_{2},...,n_{k})
	&= 
	O\left( \frac{N_{q}(J_{X})}{(q-1)}
	q^{(r^2+n_{1}^2-rn_{1}-\frac{r}{c}
		+\frac{n_{1}}{c}-1)(g-1)}\mathcal{A} \right)\\
	&+ 
	O\left(
	 S(n_{1})N_{q}(J_{X})q^{(r^{2}+n_{1}^{2}-rn_{1}
		-\frac{r}{c})(g-1)}
	\frac{\mathcal{A}}{\zeta_{X}(n_{1})}
	\right),  
\end{align*}
where 
\[
\mathcal{A}:= \zeta_{X}(2)\zeta_{X}(3)...\zeta_{X}(n_{1})\zeta_{X}(2)\zeta_{X}(3)...\zeta_{X}(r-n_{1}-1).
\]
Using Propositions \ref{log-zeta5.3} and \ref{log_Jac5.1}, we get
\begin{align*}
&
\frac{C_{L} (n_{1},n_{2},...,n_{k})}{q^{(r^2
		-\frac{r}{c})(g-1)}
	\zeta_{X}(2)\zeta_{X}(3)...\zeta_{X}(r-1)}\\
&=
 O\left(
 q^{-1} q^{(n_{1}^2
 	+\frac{n_{1}}{c}-rn_{1})(g-1)}(\log g)^{N
 	\left( 
 	\frac{2c'(r-n_{1}-1)}{q^2}+3
 	\right)} 
  \exp\left(
  \frac{2c'N}{\sqrt{q}}  (r-n_{1}-1) 
  \right) \right)\\
&+
 O\left(
 q^{(n_{1}^2-rn_{1}+1)(g-1)}(\log g)^{N\left( \frac{2c'(r-n_{1}-1)}{q^2}+3\right)} 
 \exp\left(\frac{2c'N}{\sqrt{q}} 
 \left( r-n_{1}-1\right) \right) \right).
\end{align*}
This completes the proof as the right hand side tends to $0$ for $g$ tending to $\infty.$
\end{proof}

Furthermore, using \eqref{Palpha} and Proposition \ref{Pdom5.1} we can rewrite the equation \eqref{Sigen2.2} as follows, 

\begin{rmk}\label{Siegel3}
	\begin{align*}
N_{q}(M_{L}(r,d)) 
&=
q^{(r^2-1)(g-1)}\zeta_{X}(2)\zeta_{X}(3)...\zeta_{X}(r)\\
&o\left((q-1)S(r)q^{(r^2-\frac{r}{c})(g-1)}\zeta_{X}(2)\zeta_{X}(3)...\zeta_{X}(r-1) \right).  
     \end{align*}
\end{rmk} 
\noindent
for any constant $c>1.$\\
\noindent
\textbf{Final step of the proof of Theorem \ref{thm1.3}:}\\
\noindent
From Remark \eqref {Siegel3}, 
\noindent
Let \[T_{1}:=q^{(r^2-1)(g-1)}\zeta_{X}(2)\zeta_{X}(3)...\zeta_{X}(r)\]
\noindent
and 
\[T_{2}:=o\left((q-1)S(r)q^{(r^2-\frac{r}{c})(g-1)}\zeta_{X}(2)\zeta_{X}(3)...\zeta_{X}(r-1) \right).\]
\noindent 
We choose the constant $c$ such as $1<c<r.$
Taking logarithm on both sides of $T_{1},$ we get  
\begin{equation*}
\log{T_{1}} = (r^2-1)(g-1)\log{q} + \sum\limits_{k=2}\limits^{r}\log\zeta_{X}(k).
\end{equation*}
Using Proposition \ref{log-zeta5.2} with some constant $c^{\prime}=\max\left\lbrace c_{1}, c_{2}\right\rbrace ,$ and Proposition \ref{log-zeta5.3}, we observe that
\begin{align*}
& |\log{\left(N_{q}(M_{L}(r,d))\right)}- (r^2-1)(g-1)\log q|\\
&\le c'(r-1)N\left( \frac{1}{\sqrt{q}}+ \frac{\log{\log{g}}}{q^2}\right) 
+ O\left((q-1)S(r)q^{(1-\frac{r}{c})(g-1)}\left(\log{g} \right)^{\frac{c'N}{q^{r}}}\exp\left( \frac{c'N}{\sqrt{q}}\right)\right). 
\end{align*}
Now, if we consider $A:=N\left( \frac{1}{\sqrt{q}}+ \frac{\log{\log{g}}}{q^2}\right) $
and
$B:=(\log g)^{\frac{N}{q^{2}}}\exp\left( \frac{N}{\sqrt{q}}\right),$ then $\log B=A.$
Therefore, 
\begin{equation*}
 \log{\left(N_{q}(M_{L}(r,d))\right)}- (r^2-1)(g-1)\log q \,=\, O_{r,N}(A+q^{-\sigma g}B),
 \end{equation*}
 for a suitable absolute constant $\sigma >0,$ depending on $r.$
The theorem follows upon simplifying.

\subsection{Proof of Theorem \ref{thm1.4}}  
Next we focus our attention on the family of hyperelliptic curves $\H_{\gamma,q}.$ 
Let $H$ be a hyperelliptic curve of genus $g\geq 2$ given by the affine model  $H: y^{2} = F(x)$ 
with $F$ in $\H_{\gamma,q}$. 
Suppose $(r, d)=1.$ Corresponding to a hyperelliptic curve $H,$ we use the notation $M_{L_{H}}(r,d)$ for the moduli space of stable vector bundles of rank $r$ with fixed determinant $L_{H}$ of degree $d$ defined over $\F_{q}.$
The function field $\F_{q}(H)$ corresponding to the hyperelliptic curve $H$ is a Galois extension of the rational function field $\F_{q}(x)$ of degree two.
We denote $\F_{q}(H)$ by $K'$ and $\F_{q}(x)$ by $K$ for simplicity. Let $\chi=(\frac{.}{F})$ 
denote the Legendre symbol generating 
$Gal(K'/K).$
As discussed in subsection \ref{zeta function basics}, we have
\begin{equation}\label{L-chi}
L(s,\chi)=\prod\limits_{l=1}\limits^{2g}(1-\sqrt{q}e(\theta_{l,H})q^{-s}).
\end{equation}
The Euler product of $L-$function is given by
\begin{equation}\label{Euler}
L(s,\chi) = \prod\limits_{P}(1 - \chi{(P)}\mid P \mid ^{-s})^{-1}.
\end{equation}
Taking logarithmic derivatives of (\ref{L-chi}) and (\ref{Euler}) and equating coefficients of $q^{-ms}$ for any positive integer $m,$
we get

\begin{eqnarray*}
	\sum\limits_{l=1}\limits^{2g}-e(m\theta_{l,H})
	&=& \sum\limits_{f \in \mathcal{F}_{m}}q^{-m/2}\Lambda(f)\left(\frac{f}{F}\right)
\end{eqnarray*}
where $\mathcal F_m$ is the set of all monic polynomials of degree $m$. 
Now for any $F \;\text{in}\; \hh,$ 
\begin{equation*}
\left(\frac{F}{P_\infty}\right) =
\begin{cases}
& 1  \text{ if } \deg(F)\equiv 0 \,(\text{mod}\, 2)\\
& 0   \text{ otherwise}. 
\end{cases}
\end{equation*}
Using quadratic reciprocity, we note that
\[\sum\limits_{\substack{f =\infty \\ }}q^{-m/2}\Lambda(f)\left(\frac{F} {f}\right) 
=\, q^{-m/2}\delta_{\gamma/2}. 
\]
Therefore,
\begin{equation}\label{e1}
\sum\limits_{l=1}\limits^{2g}-e(m\theta_{l,H})
= \sum\limits_{\substack{f \neq \infty  \\ \deg f=m}}q^{-m/2}\Lambda(f)\left(\frac{f}  {F}\right) 
+ q^{-m/2}\delta_{\gamma/2}.
\end{equation}
Proceeding as in the proof of Theorem \ref{thm1.3} we get 
\begin{align}\nonumber
\log{\left(N_{q}\left( M_{L_{H}}(r,d)\right) \right)} 
&-(r^2-1)(g-1)\log{q}
=
\sum\limits_{k=2}^{r}\log{\zeta_{H}(k)} \\
\label{R5.2}
&+
 O\left( 
 (q-1)S(r)q^{(1-\frac{r}{c})(g-1)}
 \left(
 \log{g} 
 \right)^{\frac{c}{q^{r}}}
 \exp\left( 
 \frac{c}{\sqrt{q}}
 \right)  
  \right)
\end{align}
for some absolute constant $1<c<r.$
For a fixed positive integer $Z$, we write
\begin{equation}\label{total sum}
\sum\limits_{k=2}^{r}\log{\zeta_{H}(k)} - \sum\limits_{k=2}^{r}(2k-1)\log{q} + \sum\limits_{k=2}^{r}\log\{(q^k-1)(q^{k-1}-1)\}=\epsilon_{1,Z} + \epsilon_{2,Z}
\end{equation}
	where,
\begin{equation}\label{gen5.1}
\epsilon_{1,Z} = 
-\sum\limits_{m\leq Z}
\left( \sum\limits_{k=2}^{r}q^{-\frac{(2k-1)}{2}m}\right)
 m^{-1}\sum\limits_{l=1}\limits^{2g}e(m\theta_{l,X}),
\end{equation}

and
\begin{equation}\label{gen5.2}
\epsilon_{2,Z} = 
-\sum\limits_{m> Z}
\left( \sum\limits_{k=2}^{r}q^{-\frac{(2k-1)}{2}m}\right) 
m^{-1}\sum\limits_{l=1}\limits^{2g}e(m\theta_{l,X}).
\end{equation}


Using \eqref{e1} in \eqref{gen5.1}, we get
\begin{eqnarray*}
	\epsilon_{1,Z} 
	&=&
	\bigtriangleup_{Z}(F) +
	\sum\limits_{m\leq Z}
	\left( \sum\limits_{k=2}^{r}q^{-km}\right)
	m^{-1}\delta_{\gamma/2},
	\end{eqnarray*}
where
 \begin{equation}\label{main5.1}
\bigtriangleup_{Z}(F):= 
\sum\limits_{m\leq Z}
\left( \sum\limits_{k=2}^{r}q^{-km}\right)
m^{-1}
\sum\limits_{\substack{f \neq \infty \\ \deg f=m}}
\Lambda(f)\left(\frac{F}{f}\right).
\end{equation}

We see,
\begin{equation}\label{upper-bounds5.1}
\sum\limits_{m\leq Z}q^{-2m}m^{-1}\sum\limits_{\substack{f \neq \infty  \\ \deg f=m}}
\Lambda(f)\left(\frac{F}{f}\right) \leq \frac{1}{q}(1+ \log{Z}).
\end{equation}
Also,
 \begin{equation*}
\sum\limits_{m\leq Z}
\bigg( 
\sum\limits_{k=2}^{r}
q^{-km}
\bigg)
m^{-1}\delta_{\gamma/2}
= -\delta_{\gamma/2}
\sum\limits_{k=2}^{r}
\log(1-1/q^k)
-\delta_{\gamma/2}
\sum\limits_{m>Z}
\frac{1}{m}
\sum\limits_{k=2}^{r}
\frac{1}{q^{km}}.
\end{equation*}
After simplification we can write, 
\begin{equation*}
\epsilon_{1,Z}= \bigtriangleup_{Z}(F) 
-\delta_{\gamma/2}
\sum\limits_{k=2}^{n}\log(1-1/q^k)
-\delta_{\gamma/2}
\sum\limits_{m>Z}\frac{1}{m}
\sum\limits_{k=2}^{r}\frac{1}{q^{km}}.
\end{equation*}
Therefore rearranging the terms in \eqref{total sum} and based on above estimates for $\epsilon_{1, Z}$, we obtain
\begin{align}
&\sum\limits_{k=2}^{r}\log{\zeta_{H}(k)} \nonumber
- \sum\limits_{k=2}^{r}(2k-1)\log{q} 
+ \sum\limits_{k=2}^{r}\log\{(q^k-1)(q^{k-1}-1)\} 
+ \delta_{\gamma/2}\log\left\lbrace\prod\limits_{k=2}^{r} (1-1/q^k)\right\rbrace\\
&=\bigtriangleup_{Z}(F) + \epsilon_{Z,F}
\label{gen5.3}
\end{align}

where, 
\begin{equation*}
\epsilon_{Z,F}=\epsilon_{2,Z} - \delta_{\gamma/2}\sum\limits_{m>Z}\frac{1}{m}\sum\limits_{k=2}^{r}\frac{1}{q^{km}}.
\end{equation*}
Using \eqref{gen5.3}, we write equation \eqref{R5.2} in a simpler form 
\begin{equation}\label{N_F5.1}
\mathcal{R}_{(r,d)}(H)-C_{q}(r) = \bigtriangleup_{Z}(F) + \epsilon_{Z}(F)
\end{equation}
where $\mathcal{R}_{(r,d)}(H),$ and $C_{q}(r)$ are as defined in \eqref{def:random variable}, and \eqref{def:log constant} respectively,
and
\begin{equation}\label{error5.1}
\epsilon_{Z}(F):= \epsilon_{Z,F} 
+ O\left(q^{(1-\frac{r}{c})(g-1)}\left(\log{g} \right)^{\frac{c}{q^{r}}}
\exp\left( \frac{c}{\sqrt{q}}\right)\right).
\end{equation}
It is easy to see that, 
\begin{equation}\label{upper-bounds5.2}
\mid\epsilon_{Z,F}\mid = \bigcirc\left(\frac{g}{Z} q^ {-3Z/2}\right).
\end{equation}

%

Choose $Z=\left[ \frac{\gamma}{3}\right].$ From \eqref{main5.1}, we write
\begin{equation}
\bigtriangleup_{\left[ \frac{\gamma}{3}\right]}(F):= \sum\limits_{k=1}^{r-1}
\R_{(\gamma, q)}^{(k)}(F)
\end{equation}

where,
\begin{equation}\label{def: random variable sequence}
\R_{(\gamma, q)}^{(k)}(F):=
\sum\limits_{m\leq \left[ \frac{\gamma}{3}\right]}
q^{-(k+1)m}m^{-1}
\sum\limits_{\substack{f \neq \infty \\ \deg f=m}}
\Lambda(f)\left(\frac{f} {F}\right)
\end{equation}
for each $1\leq k \leq r-1.$ Using the relation $g=\left[ \frac{\gamma-1}{2}\right] ,$ we conclude that
\[
\mathcal{R}_{(r,d)}(H)-C_{q}(r)
=
\sum\limits_{k=1}^{r-1}
\R_{(\gamma, q)}^{(k)}(F)
+
O(q^{-c'g})
\]
for some absolute constant $c'>0.$\\

{\bf Computation of moments and the distribution function:}\\

For any function $\psi: \H_{\gamma, q} \rightarrow \C$, we denote the mean value of $\psi$ by 
\[
\E(\psi) := 
\frac{1}{|\H_{\gamma, q}|} \sum\limits_{F\in \H_{\gamma, q}} \psi(F).
\]

Following the similar procedure as in \cite[Theorem $1.2$]{ASA}, for a positive integer $n \leq \log \gamma,$ 
one can show that for each $1\leq k \leq r-1,$  the $n$th moment of the random variable $\R_{(\gamma, q)}^{(k)}$ is given by
\begin{equation}\label{Delta-H1}
\E
\left( 
\left( 
\R_{(\gamma, q)}^{(k)}
\right) 
^n
\right) 
=
 H^{(k)}(n)+T 
\end{equation}
where, 
\begin{equation}\label{H(r)5.1}
H^{(k)}(n) := 
\sum\limits_{\substack{m_i\geq 1\\ 1\leq i\leq n}} 
\prod\limits_{i=1}\limits^{n}q^{-(k+1)m_{i}}m_{i}^{-1}
\sum\limits_{\substack{\deg f_{i}=m_{i}\\ 1\leq i \leq n\\f_{1}f_{2}...f_{n}=h^{2}}}
\Lambda(f_{1})\Lambda(f_{2})....\Lambda(f_{n})
\prod\limits_{P\mid h} (1+ \mid P \mid ^{-1})^{-1},
\end{equation}
and $T=O(q^{-c''\gamma})$ for some $c''>0.$
Therefore, if $q$ is fixed then for each fixed positive integer $n$,
\begin{equation}\label{weak convergence}
\lim\limits_{\gamma\rightarrow \infty} 
\E
\left(
\left( 
\R_{(\gamma, q)}^{(k)}
\right) 
^n
\right)
 = H^{(k)}(n),
\end{equation}
for each $1\leq k\leq r-1.$

Now proceeding similarly as in \cite[Proposition $1$, $3$]{XiZa} one can prove the following two results. For the sake of completeness we give a brief outline of the proof here.
\begin{prop}\label{H(r)1}
	For any positive integer $n \geq 1$, we have
	$$ H^{(k)}(n)= \sum_{s=1}^{n} \frac{n!}{2^{s}s!}
	\sum_{\substack{\sum\limits_{i=1}^{s}\lambda_{i}=n \\ \lambda_{i}\geq 1}} 
	\sum_{\substack{P_{1},...,P_{s}\\\text{distinct} }}
	\prod_{i=1}^{s}\frac{u_{P_{i}}^{\lambda_{i}}+ (-1)^{\lambda_{i}}v_{P_{i}}^{\lambda_{i}}}
	{\lambda_{i}!(1+ \mid P_{i}\mid ^{-1})}$$
	where the sum on the right hand side is over all positive integer 
	$\lambda_{i}, i= 1, 2, ..., s$ such that $\sum\limits_{i=1}\limits^{s}\lambda_{i}=n,$ 
	and over all distinct monic, irreducible polynomials $P_{i} \;\text{in}\; \F_{q}[x]$ with 
	\begin{eqnarray*}
		u_{P_{i}} &=& -\log (1-|P_{i}|^{-(k+1)}),\\
		v_{P_{i}} &=& \log(1+ |P_{i}|^{-(k+1)}).
	\end{eqnarray*}
\end{prop}
\begin{proof}
From \eqref{H(r)5.1}, we can write
\begin{equation*}
H^{(k)}(n)=\sum\limits_{h}
\prod\limits_{P|h}
(1+|P|^{-1})^{-1}
|h|^{-2(k+1)}
\sum\limits_{\substack{\deg f_{i}=m_{i}\\ 1\leq i \leq n\\f_{1}f_{2}...f_{n}=h^{2}}}
\prod\limits_{i=1}^{n}
\frac{\Lambda(f_{i})}
{\deg f_{i}}.
\end{equation*}
Now, the inside sum survives only when each $f_{i}$ is a power of prime, and hence $\omega(h)\leq n,$ where $\omega(h)$ denotes the number of distinct prime factors of $h.$ Therefore, partitioning the above sum depending on the number of distinct prime factors in $h,$ we write
\begin{equation}\label{def: H(s,n)}
H^{(k)}(n)=\sum\limits_{s=1}^{n}
H^{(k)}(s,n),
\end{equation}  
where
\begin{equation*}
H^{(k)}(s,n)
=
\sum\limits_{\substack{h\\ \omega(h)=s}}
\prod\limits_{P|h}
(1+|P|^{-1})^{-1}
|h|^{-2(k+1)}
\sum\limits_{\substack{\deg f_{i}=m_{i}\\ 1\leq i \leq n\\f_{1}f_{2}...f_{n}=h^{2}}}
\prod\limits_{i=1}^{n}
\frac{\Lambda(f_{i})}
{\deg f_{i}}.
\end{equation*}
For each class $H^{(k)}(s,n),$ we can choose a tuple $(P_{1}, P_{2},..,P_{s})$ of distinct primes and their corresponding exponents $(\alpha_{1}, \alpha_{2},...,\alpha_{s}).$ Also, the ordering of each $P_{i}$ is irrelevant for a fix $h.$ Therefore, we obtain
\begin{equation*}
H^{(k)}(s,n)
=
\frac{1}{s!}
\sum\limits_{\substack{P_{1},...,P_{s}\\ \text{distinct}      }}
\sum\limits_{\substack{\alpha_{i} \geq 1\\ 1\leq i\leq s   \\
h=\prod P_{i}^{\alpha_{i}}   }}
\prod\limits_{i=1}^{s}
(1+|P_{i}|^{-1})^{-1}
|P_{i}|^{-2(k+1)\alpha_{i}}
\sum\limits_{\substack{\deg f_{i}=m_{i}\\ 1\leq i \leq n\\f_{1}f_{2}...f_{n}=h^{2}}}
\prod\limits_{j=1}^{n}
\frac{\Lambda(f_{j})}
{\deg f_{j}}.
\end{equation*}

Now, for $1\leq j\leq n,$ each $f_{j}=Q_{i}^{\beta_{i}}$ for some prime $Q_{i} \in \left\lbrace P_{1}, P_{2},...,P_{s}\right\rbrace $ such that $s\leq n,$ and $\beta_{i}\geq 1,$ and
\[
\sum\limits_{j\in A_{i}}
\beta_{j}
=
2\alpha_{i}
\]
where $A_{i}$ is the set of all indices $j$ such that $f_{j}$ is a power of $P_{i}.$ Therefore, for a fixed set of distinct prime $(P_{1}, P_{2},..,P_{s})$, and the corresponding exponent $(\alpha_{1}, \alpha_{2},...,\alpha_{s})$, we have
\begin{equation*}
\prod\limits_{j=1}^{n}
\frac{\Lambda(f_{j})}
{\deg f_{j}}
=
\prod\limits_{j=1}^{n}
\frac{1}
{\beta_{j}}.
\end{equation*}
Suppose $|A_{i}|=\lambda_{i} \geq 1$ for $1\leq i\leq s.$ Therefore, $\sum\limits_{i=1}^{s}\lambda_{i}=n.$ Hence we can write
\begin{equation*}
H^{(k)}(s,n)
=
\frac{1}{s!}
\sum\limits_{\substack{ \sum\limits_{i=1}^{s}\lambda_{i}=n\\
\lambda_{i}\geq 1  }}
\frac{n!}{\lambda!}
\sum\limits_{\substack{P_{1},...,P_{s}\\ \text{distinct}      }}
\prod\limits_{i=1}^{s}
\left( 
(1+|P_{i}|^{-1})^{-1}
\sum\limits_{\substack{  \sum\limits_{j=1}^{\lambda_{i}} a_{j} \equiv 0 \pmod* 2\\
		a_{j}\geq 1    }}
	\frac{  |P_{i}|^{-(k+1)  \sum\limits_{j=1}^{\lambda_{i}} a_{j} }  }
	{\prod\limits_{j=1}^{\lambda_{i}} a_{j} }
\right). 
\end{equation*}

Following a similar argument as in \cite[Proposition 1]{XiZa}, it can be shown that
\begin{equation*}
\sum\limits_{\substack{  \sum\limits_{j=1}^{\lambda} a_{j} \equiv 0 \pmod* 2\\
		a_{j}\geq 1    }}
\frac{  |P|^{-(k+1)  \sum\limits_{j=1}^{\lambda} a_{j} }  }
{\prod\limits_{j=1}^{\lambda} a_{j} }
=
\frac{1}{2}
\left( 
u_{P}^{\lambda}
+
(-1)^{\lambda}
v_{P}^{\lambda}
\right),  
\end{equation*}
for any positive integer $\lambda \geq 1,$ and each prime $P \in \F_{q}[x],$ where
\begin{eqnarray*}
	u_{P} &=& -\log (1-|P|^{-(k+1)}),\\
	v_{P} &=& \log(1+ |P|^{-(k+1)}).
\end{eqnarray*}
Hence the result.
\end{proof}

Now suppose $\R_{k}$ is a random variable such that  for any positive integer $n,$
 \[
 \mathbb{E}
 \left( \R_{k}^n\right)=H^{(k)}(n).
 \]
 The characteristic function $\phi_{\R_{k}}(t)$ of $\R_{k}$ is given by
\[
\phi_{\R_{k}}(t)
=
\mathbb{E}
\left(
e^{it\R_{k}}
 \right), 
\]
where $t\in \mathbb{R},$ and $i$ is the imaginary unit.
Writing the characteristic function $\phi_{\R_k}(t)$ in terms of the $n^{\text{th}}$ moment $H^{(k)}(n),$ and applying Proposition \ref{H(r)1}, we see  
\begin{eqnarray*}  
	\phi_{\R_k}(t) 
	 &=& 
	1
	+
	 \sum_{n=1}^{\infty} \frac{(it)^{n}}{n!}
	 H^{(k)}(n)\\
	 &=&
	 1+
	 \sum_{n=1}^{\infty} \frac{(it)^{n}}{n!}
	 \sum_{s=1}^{n} \frac{n!}{2^{s}s!}
	\sum_{\substack{\sum\limits_{j=1}^{s}\lambda_{j}=n\\ \lambda_{j}\geq 1}} 
	\sum_{\substack{P_{1},...,P_{s}\\ \text{distinct}      }}
	\prod_{j=1}^{s}\frac{u_{P_{j}}^{\lambda_{j}}
		+ (-1)^{\lambda_{j}}v_{P_{j}}^{\lambda_{j}}}{\lambda_{j}!(1+ \mid P_{j}\mid ^{-1})}.
\end{eqnarray*}
 
Changing the order of summation we get,
\begin{eqnarray*}
	\phi_{\R_k}(t)
	&=&
	 1
	+
	 \sum\limits_{s=1}\limits^{\infty}\frac{1}{2^{s}s!}
	\sum\limits_{\substack{P_{1},...,P_{s}\\ \text{distinct}      }}
	\prod\limits_{j=1}\limits^{s}
	\left(
	\sum\limits_{\lambda_{j}=1}\limits^{\infty}
	\frac{(it)^{\lambda_{j}}
		\left(
		u_{P_{j}}^{\lambda_{j}}+ (-1)^{\lambda_{j}}v_{P_{j}}^{\lambda_{j}}
		\right)}
	{\lambda_{j}!(1+ \mid P_{j}\mid ^{-1})}
	\right)\\
	&=& 
	1+
	 \sum\limits_{s=1}\limits^{\infty}
	 \frac{1}{2^{s}s!}
	 \sum\limits_{\substack{P_{1},...,P_{s}\\ \text{distinct}      }}
	 \prod\limits_{j=1}\limits^{s}
	 \left(
	 \frac{(1-\mid P_{j} \mid^{-(k+1)})^{-it}
	 	+ 
	 	(1+\mid P_{j} \mid^{-(k+1)})^{-it}-2}{(1+ \mid P_{j}\mid ^{-1})}
	 \right).
\end{eqnarray*}

This completes the proof of $(1)$ of Theorem \ref{thm1.4}.

The next result gives an asymptotic formula for $H^{(k)}(n)$ for each $1\leq k\leq r-1,$ when $q\rightarrow \infty.$
\begin{prop}\label{H(r)2}
As $q\rightarrow \infty,$ the $n$-th moment
\[
H^{(k)}(n) = \frac{\delta_{n/2}n!}{2^{n/2}(n/2)!}q^{\frac{-(2k+1)n}{2}} 
+ \bigcirc_{n}
\left( 
q^{ \frac{-(2k+1)n}{2}-1 }  
 \right) .
\]
\end{prop}
\begin{proof}
	This is an exact analogue \cite[Proposition 3]{XiZa} and we skip the proof. 
\end{proof}
Considering 
$q^{\frac{(2k+1)}{2}}\R_{k}$ as a random variable on the space $\hh$, as both $\gamma,$ and $q \rightarrow \infty,$ we see that
all its moments are asymptotic to the corresponding moments of a 
standard Gaussian distribution where the odd moments vanish and the even moments are 
\begin{equation*}
\frac{1}{\sqrt{2\pi}}\int_{-\infty}^{\infty}\tau^{2n}e^{-\tau^{2}/2} d\tau = \frac{(2n)!}{2^{n}n!}.
\end{equation*} 
Hence the corresponding characteristic function converges to characteristic function of Gaussian distribution. Finally using Continuity theorem ( see Theorem $3.3.6$ in \cite{Du}), we obtain result $(3)$ of Theorem \ref{thm1.4}. \\

\noindent
\textbf{Computation of Co-variance:}

In this section we prove that the random variables $\R_{(\gamma, q)}^{(k)}$ over $\hh$ for $1\leq k \leq r-1,$ are not independent.
%
Before going to the details of the proof, first we would list a set of results required in the proof.
\begin{lem}\label{character-sum}
	Let $h$ be a polynomial in $\F_{q}[x]. $ For any non-trivial Dirichlet character $\chi\pmod h $, 
	we have 
	\[
	\frac{1}{\#\H_{\gamma,q}}\sum\limits_{F \in \H_{\gamma,q}}\chi(F) 
	\leq\, \frac{2^{\deg h }-1}{(1-1/q)q^{\gamma/2}}.
	\]
\end{lem}
\begin{proof}
	The proof follows from [Lemma 3.1 of \cite{FaRu}].
\end{proof}
\begin{lem}\label{trivial-sum}
	Let $h \;\text{in}\; \F_{q}[x]$ be a monic square-free polynomial, then 
	\begin{eqnarray*}
		\frac{1}{\#\H_{\gamma,q}}\sum \limits_{\substack{F\in \h\\gcd(F,h)=1}}1
		&=& 
		\prod\limits_{p\mid h} (1+ \mid p \mid^{-1})^{-1} 
		+ \bigcirc(q^{-\gamma/2} \tau(h))
	\end{eqnarray*}
	where $\tau(h) = \sum\limits_{D\mid h}1.$
\end{lem}

\begin{proof}
	This is essentially [\cite{Ru},Lemma 5]. For more details one can see [\cite{XiZa}, Lemma 2].
\end{proof}    
Using above mentioned results we are ready to compute the covariance of the random variables. The following proposition proves the statement of the Theorem\ref{thm1.4}(2). 
\begin{prop}\label{thm covariance}
	The limiting covariance of the random variables $\R_{(\gamma, q)}^{(i)},\, R_{(\gamma, q)}^{(j)}$ for $1\leq i \neq j \leq r-1$ is
	\[
	\lim\limits_{\gamma \rightarrow \infty}
	Cov
	\left( 
	\R_{(\gamma, q)}^{(i)},\, R_{(\gamma, q)}^{(j)}
	\right) 
	=
	q^{-(i+j+1)} + O(q^{-(i+j+2)}).
	\]
\end{prop}

\begin{proof}
	
	We know that 
	\begin{align*}
	Cov
	\left( 
	\R_{(\gamma, q)}^{(i)},\, R_{(\gamma, q)}^{(j)}
	\right) 
	&=
	\E
	\left( 
	\R_{(\gamma, q)}^{(i)}
	\R_{(\gamma, q)}^{(j)}
	\right) 
	-
	\E
	\left( 
	\R_{(\gamma, q)}^{(i)}
	\right) 
	\E
	\left( 
	\R_{(\gamma, q)}^{(j)}
	\right) .
	\end{align*}
	We have
	\begin{align*}
	\E
	\left( 
	\R_{(\gamma, q)}^{(i)}
	\R_{(\gamma, q)}^{(j)}
	\right) 
	&=
	\sum\limits_{n,m \leq Z}
	\frac{q^{-n(i+1)-m(j+1)}}{nm}
	\sum\limits_{\substack{\deg f = n \\ \deg g=m}}
	\Lambda(f)\Lambda(g)
	\E
	\left(  
	\left(\frac{fg}{\cdot}\right)
	\right),
	\end{align*}
	where $Z=\left[ \frac{\gamma}{3}\right].$ We first consider the case that $fg$ is not a square in $\F_{q}[x]$. 
	Then using quadratic reciprocity we see that $\left(\frac{\cdot}{fg}\right): \F_{q}[x] \rightarrow \C$ is a nontrivial Dirichlet 
	character modulo $fg.$
	Let $T_1$ be the total contribution from this case to $\E
	\left( 
	\R_{(\gamma, q)}^{(i)}
	\R_{(\gamma, q)}^{(j)}
	\right)$. Using lemma \ref{character-sum} we obtain
	\begin{align*}
	\quad T_{1}
	&\leq
	\sum\limits_{n, m \leq Z} 
	\frac{q^{-n(i+1)-m(j+1)}}{nm}
	\sum\limits_{\substack{\deg f = n \\ \deg g=m}}
	\Lambda(f)\Lambda(g)
	\frac{2^{m+n-1}}
	{(1-1/q)q^{\gamma/2}}\\
	&\leq 
	\sum\limits_{n, m \leq Z} 
	\frac{q^{-in-jm}}  {nm}
	\frac{2^{2Z}q^{-\gamma/2}}
	{2(1-1/q)}
	\end{align*}
	where the last inequality we get using 
	\[
	\sum\limits_{\substack{\deg f=n\\f\neq \infty}}
	\Lambda(f)
	=q^{n}.
	\]
	Finally, since $i, j \geq 1,$
	$T_{1}$ can be bounded above as
	\begin{equation}
	T_{1} \leq 
	\frac{1}{(q-1)^2}\frac{2^{2Z}q^{-\gamma/2}}{2(1-1/q)}
	\leq
	q^{-\gamma/2-2}2^{2Z}. 
	\end{equation}
	Next we suppose that  $fg$ is a square in $\F_{q}[x]$ and 
	$fg= h^{2}$. We write the square free part of $h$ as  
	$ \tilde{h}=\prod\limits_{\substack{P\mid h\\ P\,\text{ prime}}}P.$ 
	Using lemma \ref{trivial-sum} we get
	\begin{eqnarray*}
		\E
		\left( 
		\left( \frac{\cdot}{h^2}\right)  
		\right) 
		= \frac{1}{\#\H_{\gamma,q}}\sum \limits_{\substack{F\in \h\\gcd(F,\tilde{h})=1}}1
		&=& \prod\limits_{p\mid \tilde{h}} (1+ \mid p \mid^{-1})^{-1} + \bigcirc(q^{-\gamma/2} \tau(\tilde{h})).
	\end{eqnarray*}
	Let $T_2$ be the total contribution from the error term $\bigcirc(q^{-\gamma/2} \tau(\tilde{h}))$ to $\E
	\left( 
	\R_{(\gamma, q)}^{(i)}
	\R_{(\gamma, q)}^{(j)}
	\right) $. Then
	\begin{align*}
	\quad T_{2}
	&=
	\sum\limits_{n, m \leq Z} 
	\frac{q^{-n(i+1)-m(j+1)}}{nm}
	\sum\limits_{\substack{\deg f = n \\ \deg g=m}}
	\Lambda(f)\Lambda(g)
	q^{-\gamma/2}\tau{(\tilde{h})}.
	\end{align*}
	Since $f,$ and $g$ are prime powers in the second sum, $h$ and therefore $\tilde{h}$ has at most $2$ distinct prime factors.
	Therefore 
	$\tau(\tilde{h}) \leq 2$. We see that
	\begin{align*}
	T_{2}
	&\leq
	q^{-\gamma/2}
	\left( 
	\sum\limits_{n\leq Z} 
	\frac{q^{-n(i+1)}}{n}
	\sum\limits_{\substack{\deg f = n}}
	\Lambda(f)
	\right) 
	\left( 
	\sum\limits_{n\leq Z} 
	\frac{q^{-m(j+1)}}{m}
	\sum\limits_{\substack{\deg g = m}}
	\Lambda(g)
	\right) \\
	&\leq
	q^{-\gamma/2}
	\left( 
	\sum\limits_{n\leq Z}
	\frac{q^{-in}}{n}
	\right)
	\left( 
	\sum\limits_{m\leq Z}
	\frac{q^{-jm}}{m}
	\right)\\
	&\leq 
	q^{-\gamma/2}
	\frac{1}{(q-1)^2}\\
	&\leq
	q^{-\gamma/2-2}.
	\end{align*}
	\noindent
	Let $M$ be the total contribution from the main term 
	$\prod\limits_{p\mid \tilde{h}} (1+ \mid p \mid^{-1})^{-1}$ 
	to $\E
	\left( 
	\R_{(\gamma, q)}^{(i)}
	\R_{(\gamma, q)}^{(j)}
	\right) $. We estimate 
	\begin{eqnarray*}
		\quad M
		=\sum\limits_{n, m \leq Z} 
		\frac{q^{-n(i+1)-m(j+1)}}{nm}
		\sum\limits_{\substack{fg=h^2\\ \deg f = n \\ \deg g=m}}
		\Lambda(f)\Lambda(g)
		\prod\limits_{p\mid h} (1+ \mid p \mid^{-1})^{-1}.
	\end{eqnarray*}
	Next, we remove the dependence on $Z$ in the first sum and extend it to all $n, m\geq 1.$ Therefore, 
	the condition $fg=h^2$ remains,
	causing an error bounded by 
	
	\begin{equation*}
	T_{3}
	\,:=\,
	\sum\limits_{\substack{h\\ \deg h>Z/2}}
	\prod\limits_{p\mid h} (1+ \mid p \mid^{-1})^{-1}
	\sum\limits_{\substack{f, g\\fg= h^{2}}}
	\frac{\Lambda(f)\Lambda(g)}
	{\deg f \deg g}
	|f|^{-(i+1)}  |g|^{-(j+1)}
	\end{equation*}

\noindent
To find upper bound for $T_3$, we note that $\Lambda(f)\le \deg(f)$. Moreover the $f$ and $g$ appearing in $T_3$
are all powers of primes. Hence the outer sum is over monic polynomials $h$ having atmost $2$ distinct
prime factors. Therefore given an $h$ the number of choices for each $f$ is atmost $2\deg(h)$. Considering $i<j,$ we get
\begin{align*}
T_{3}
&\leq
\sum\limits_{\substack{h\\ \deg h>Z/2}}
2\mid h\mid^{-2(i+1)}
2 (2\deg(h))^2\\
&\ll 
\sum\limits_{k> Z/2}
\sum\limits_{\deg h =k}
|h|^{-2(i+1)}
(\deg h)^2\\
&\leq
\sum\limits_{k> Z/2}
k^2 q^{-2k(i+1)} q^k\\
&\leq
q^{-5Z/4}
\end{align*}
using $i\geq 1,$ and $k^2<q^{k/2}$ for $q$ large enough.
Combining the above estimates together we obtain 
\begin{equation}\label{Delta-H}
\E
\left( 
\R_{(\gamma, q)}^{(i)}
\R_{(\gamma, q)}^{(j)}
\right) 
=M_{1} + T
\end{equation}
where 
\begin{equation}
M_{1}
=
\sum\limits_{h}
\prod\limits_{P|h}
(1+|P|^{-1})^{-1}
\sum\limits_{\substack{f,g\\fg=h^2}}
\frac{\Lambda(f)\Lambda(g)}
{\deg f \deg g}
|f|^{-(i+1)}|g|^{-(j+1)}
\end{equation}
and
$T=T_{1} + T_{2} +  T_{3}= O(q^{-Z})=O(q^{-\gamma/3}).$ 

Next, since number of prime factor of $h$ could be atmost $2$ to survive the inside sum, depending on the number of prime factors of $h,$ we partition the sum
\[
M_{1}=M_{1,1} +M_{1,2}
\]
where
\begin{align*}
M_{1,1}
&:=
\sum\limits_{P}
\sum\limits_{\alpha \geq 1}
\sum\limits_{\substack{\alpha_{1}, \alpha_{2} \\ \alpha_{1}+\alpha_{2}=2\alpha  }}
\frac{
	\Lambda(P^{\alpha_{1}})\Lambda(P^{\alpha_{2}})
}
{\alpha_{1}\alpha_{2}(\deg P)^2}
\frac{
	|P|^{-(i+1)\alpha_{1} -(j+1)\alpha_{2}}
}
{
	(1+|P|^{-1})
}\\
&=
\sum\limits_{P}
(1+|P|^{-1})^{-1}
\sum\limits_{\substack{ \alpha_{1} +\alpha_{2} \equiv 0 \pmod* 2 \\ \alpha_{i} \geq 1     }}
\frac{|P|^{-(i+1)\alpha_{1} -(j+1)\alpha_{2}   }}   {  \alpha_{1} \alpha_{2}  }.
\end{align*}
and
\[
M_{1,2}:=
\sum\limits_{P_{1}\neq P_{2}}
\sum\limits_{\substack{\alpha_{1} \equiv 0 \pmod* 2\\ \alpha_{2} \equiv 0 \pmod* 2 \\ \alpha_{i} \geq 1   } }
\frac{
	\Lambda(P_{1}^{\alpha_{1}})\Lambda(P_{2}^{\alpha_{2}})
}
{\alpha_{1}\alpha_{2}(\deg P_{1}) (\deg P_{2})}
\frac{
	|P_{1}|^{-(i+1)\alpha_{1} }|P_{2}|^{-(j+1)\alpha_{2}}
}
{
	(1+|P_{1}|^{-1})
	(1+|P_{2}|^{-1})
}.
\]
We define for $1\leq k\leq r-1,$
\begin{equation*}
\eta_{P^{k}}:=
\sum\limits_{\substack{ \alpha \equiv 0 \pmod* 2 \\ \alpha \geq 1     }}
\frac{|P|^{-(k+1)\alpha   }}   {  \alpha  },
\end{equation*}
and
\begin{equation*}
\tau_{P^{k}}:=
\sum\limits_{\substack{ \alpha \equiv 1 \pmod* 2 \\ \alpha \geq 1     }}
\frac{|P|^{-(k+1)\alpha    }}   {  \alpha  }.
\end{equation*}
Therefore, we have
\begin{equation*}
M_{1,1}=\sum\limits_{P}
\frac{  \eta_{P^{i}}\eta_{P^j} + \tau_{P^{i}}\tau_{P^j}   }
{(1+|P|^{-1})},
\end{equation*} 
and
\begin{equation*}
M_{1,2}=\sum\limits_{P_{1} \neq P_{2}}
\frac{  \eta_{P_{1}^{i}}\eta_{P_{2}^{j}}  }
{(1+|P_{1}|^{-1}) (1+|P_{2}|^{-1})               }.
\end{equation*}


For a fixed $q,$ using \eqref{weak convergence} and Proposition \ref{H(r)1}, we obtain
\begin{equation*}
\lim\limits_{\gamma \rightarrow \infty}
\E
\left( 
\R_{(\gamma, q)}^{(k)}
\right) 
=
H^{(k)}(1)
=-\sum\limits_{P}
\frac{ \log \left(  1-|P|^{-2(k+1)}    \right)     }
{2 (1 + |P|^{-1})}
\end{equation*}
for all $1\leq k \leq r-1$.
Therefore,
\begin{align*}
\lim\limits_{\gamma \rightarrow \infty}
\E
\left( 
\R_{(\gamma, q)}^{(i)}
\right) 
\E
\left( 
\R_{(\gamma, q)}^{(j)}
\right) 
&=
-\sum\limits_{P_{1}, P_{2}}
\frac{ \log \left(  1-|P_{1}|^{-2(i+1)}    \right)    \log \left(  1-|P_{2}|^{-2(j+1)}    \right) }
{4 (1 + |P_{1}|^{-1}) (1 + |P_{2}|^{-1})         }\\
&=
\sum\limits_{P}
\frac{  \eta_{P^{i}}\eta_{P^{j}}   }
{(1+|P|^{-1})^2}
+
M_{1,2}.
\end{align*}

We obtain, for $1 \leq i\neq j \leq r-1,$
\begin{align}
\nonumber
\lim\limits_{\gamma \rightarrow \infty}
\E
\left( 
\R_{(\gamma, q)}^{(i)}
\R_{(\gamma, q)}^{(j)}
\right) 
-
\lim\limits_{\gamma \rightarrow \infty}
\E
\left( 
\R_{(\gamma, q)}^{(i)}
\right) 
\E
\left( 
\R_{(\gamma, q)}^{(j)}
\right) 
&=\sum\limits_{P}
\left\lbrace 
\frac{  \eta_{P^{i}}\eta_{P^{j}} + \tau_{P^{i}}\tau_{P^{j}}   }
{(1+|P|^{-1})}
-
\frac{  \eta_{P^{i}}\eta_{P^{j}}   }
{(1+|P|^{-1})^2}
\right\rbrace \\
\label{covariance1}
&=
\sum\limits_{P}
\left\lbrace 
\frac{   \tau_{P^{i}}\tau_{P^{j}}   }
{(1+|P|^{-1})}
+
\frac{  \eta_{P^{i} }\eta_{P^{j}}   }
{|P|(1+|P|^{-1})^2}
\right\rbrace.
\end{align}

Using the Taylor series expansion
\begin{equation*}
-\log(1-x)=\sum\limits_{n\geq 1}
\frac{x^n}{n},
\qquad |x|<1,
\end{equation*}
we find that for any $1\leq k \leq r-1,$
\begin{equation*}\label{recurence sum eta}
\eta_{P^{k}}=-\frac{1}{2}
\left( 
\log(1-|P|^{-2(k+1)})
\right) 
\end{equation*}
and
\begin{equation*}\label{recurence sum}
\eta_{P^{k}}+\tau_{P^{k}}
=
-
\log {(1-|P|^{-(k+1)})},
\end{equation*}
which implies
\begin{equation*}\label{recurence sum tau}
\tau_{P^{k}}=
-
\log {(1-|P|^{-(k+1)})}
+\frac{1}{2}
\left( 
\log(1-|P|^{-2(k+1)})
\right). 
\end{equation*}

We further estimate $\eta_{P^{k}}$ and $\tau_{P^{k}}$, as follows
\[
\tau_{P^{k}}=|P|^{-(k+1)} +O(|P|^{-3(k+1)}),
\]
and
\[
\eta_{P^{k}}=\frac{|P|^{-2(k+1)}}{2} +O(|P|^{-4(k+1)}),
\]
for all $k.$
For any integer $n\geq 1,$ putting these bounds in \eqref{covariance1}, we obtain

\begin{align}
\nonumber
\lim\limits_{\gamma \rightarrow \infty}
\E
\left( 
\R_{(\gamma, q)}^{(i)}
\R_{(\gamma, q)}^{(j)}
\right) 
-
\lim\limits_{\gamma \rightarrow \infty}
\E
\left( 
\R_{(\gamma, q)}^{(i)}
\right) 
\E
\left( 
\R_{(\gamma, q)}^{(j)}
\right) 
&=
\sum\limits_{P}
|P|^{-(i+j+2)} + O_{n}(|P^{-(i+j+3)}|)\\
\nonumber
&=
\sum\limits_{n\geq 1}
\sum\limits_{\deg P=n}
|P|^{-(i+j+2)} + O_{n}(|P^{-(i+j+3)}|)\\
\label{covariance2}
&=
\sum\limits_{n\geq 1}
\left( 
q^{-n(i+j+2)} + O_{n}(q^{-n(i+j+3)})
\right) 
\pi_{q}(n)
\end{align}
where
$\pi_{q}(n)$ denote the number of monic, irreducible polynomials in $\F_{q}[x]$ of degree $n\geq 1.$ Note that, from the prime number theorem for polynomials (cf. \cite[Theorem $2.2$]{Ro}), we have
\begin{equation}\label{PNT}
\pi_{q}(n)=
\frac{q^n} {n} 
+
O
\left( 
\frac{q^{n/2}} {n}
\right). 
\end{equation}
For $n\geq 2,$ using \eqref{PNT}, and using the fact $\pi_{q}(1)=q$ in \eqref{covariance2}, we obtain the proposition.

\end{proof}

\section{Distribution on $M^{s}_{\mathcal{O}_{H}}(2,0)$}\label{distribution Ms}

In this section, all the asymptotic formulas we will be considering hold for the hyperelliptic curves having sufficiently large genus.
\subsection{Proof of Theorem \ref{thm1.5}}

Over the family of hyperelliptic curves $\hh,$ we recall from Proposition \ref{Ms2.1} that, the number of $\F_{q}$-rational points over the moduli space of stable vector bundles $M^{s}_{\mathcal{O}_{H}}(2,0)$ defined over the smooth projective hyperelliptic curve $H: y^{2}=F(x)$ in $\hh,$ is given by the following expression:
\begin{equation}\label{points over Ms}
N_{q}(M^{s}_{\mathcal{O}_{H}}(2,0))= q^{3g-3}\zeta_{H}(2)
-\frac{\left( q^{g+1}-q^2+q\right) }{(q-1)^2(q+1)}N_{q}(J_{H})
-\frac{1}{2(q+1)}N_{q^2}(J_{H})
+ \frac{1}{2(q+1)}2^{2g}.
\end{equation}

For simplicity we denote,
\[
T_{1}:= q^{3g-3}\zeta_{H}(2)
\]
and
\[
T_{2}:= \frac{\left( q^{g+1}-q^2+q\right) }{(q-1)^2(q+1)}N_{q}(J_{H})
+ \frac{1}{2(q+1)}N_{q^2}(J_{H}) 
- \frac{1}{2(q+1)}2^{2g}
\]
Using Proposition \ref{log_Jac5.1} for hyperelliptic curves $(N=2)$, we see
\[T_{2}=O(q^{2g}(\log{g})^{c}),\]
for some absolute constant $c>0.$
Next using Proposition \ref{log-zeta5.3}, we obtain
\begin{equation*}
N_{q}(M^{s}_{\mathcal{O}_{H}}(2,0))
=T_{1}\left( 1 - \frac{T_{2}}{T_{1}}\right) 
=q^{3g-3}\zeta_{H}(2)\left( 1+O\left( q^{-g}(\log{g})^{c'} \right)\right) 
\end{equation*}
for some $c^{\prime}>c.$
Taking logarithm on both sides of the above equation, we have
\begin{equation}\label{eqn Ms2}
\log{N_{q}(M^{s}_{\mathcal{O}_{H}}(2,0))} - (3g-3)\log q =\log{\zeta_{H}(2)}+  O\left( q^{-g}(\log{g})^{c'}\right).
\end{equation} 
In \eqref{R5.2}, if we put $r=2,$ we observe that the right hand side of \eqref{eqn Ms2} matches with \eqref{R5.2} for a suitable choice of $c$. As a consequence, Theorem \ref{thm1.5} follows from Theorem \ref{thm1.4} for the case $r=2.$
Alternatively, For $g \rightarrow\infty,$ Theorem \ref{thm1.5} follows by invoking \cite[Theorem $1.2$]{ASA}.

\section{Distribution on $\widetilde {N}_{\mathcal{O}_{H}}(4,0)$} \label{distribution N tilde}
\subsection{Proof of Theorem \ref{thm1.6}}

We recall from equation \eqref{pointN2.2} that over a smooth projective hyperelliptic curve $H: y^{2}=F(x)$ in $\hh,$
\begin{equation}\label{pointN}
N_{q}({\NN}) \,=\, N_{q}(M^{s}_{\mathcal{O}_{H}}(2,0)) + N_{q}(Y) + 2^{2g}N_{q}(R)+ 2^{2g}N_{q}(S).
\end{equation}
Further, from \eqref{pointY} we write the $\F_{q}$-rational points on $Y$ as follows,
\begin{equation*}
N_{q}(Y) = \left(  \frac{q^{2g-3}-q}{2(q-1)(q+1)}\right)N_{q}(J_{H}) 
+\left(\frac{q^{2g-2}-1}{2(q-1)(q+1)} \right)N_{q^2}(J_{H})
-\left( \frac{q^{2g-3}-1}{2(q-1)}\right)2^{2g}. 
\end{equation*}
Also, cardinalities of the grassmannians are given by
\begin{equation*}
N_{q}(G(2,g)) = \frac{(q^{g}-1)(q^{g-1}-1)}{(q-1)^{2}(q+1)},
\end{equation*}
and,
\begin{equation*}
N_{q}(G(3,g)) = \frac{(q^{g}-1)(q^{g-1}-1)(q^{g-2}-1)}{(q^3-1)(q-1)^{2}(q+1)}.
\end{equation*}
Therefore,
\begin{align*}
	N_{q}(G(2,g))+N_{q}(G(3,g)) 
	&= q^{3g-9}\left( 1+ O\left( \frac{1}{q^{g-5}}\right) \right). 
\end{align*}
Using \eqref{points over Ms}, along with the above results in \eqref{pointN}, and further simplifying, we can write
\begin{equation}\label{eqn N}
N_{q}(\NN)\,=\, U_{1}+U_{2},
\end{equation}
where
\[
U_{1}:=\frac{1}{2}q^{2g-4}N_{q^2}(J_{H})
\left( 1+O\left( \frac{1}{q}\right) \right),
\]
and 
\[
U_{2}:=q^{3g-3}\log{\zeta_{H}(2)} + \frac{1}{2}q^{2g-5}
N_{q}(J_{H})
\left( 1+O\left( \frac{1}{q}\right) \right)  
+ q^{3g-9}
2^{2g}
\left( 1+O\left( \frac{1}{q}\right) \right).
\]
Using Proposition \ref{log-zeta5.2} and Proposition \ref{log_Jac5.1}, we get 
\begin{equation*}
U_{2}= O\left( q^{\frac{7g}{2}}\right) 
\end{equation*}
and 
\begin{equation*}
U_{1}\geq q^{4g-4}\left( \log (7g)\right)^{-3}.
\end{equation*}
For large $q$ and $g$ we have ${\left( \frac{U_{2}}{U_{1}}\right) }=O(q^{-g/2}\left( \log(7g)\right)^3 ).$

Finally from \eqref{eqn N}, we write
\begin{equation*}
N_{q}({\NN})=U_{1}
\left( 1+O(q^{-g/2}
\left( \log(7g)
\right)^3)
\right). 
\end{equation*}
Taking logarithm on both sides of the above equation and using Proposition \ref{Xiong-Zaha}, we find that
\begin{equation}\label{desingular}
\log{N_{q}({\NN})}= \log{U_{1}}+O(q^{-g/2}\left( \log(7g)\right)^3)=(4g-4)\log q + O(1).
\end{equation} 
Next, we recall the following result of Xiong and Zaharescu (Theorem 2 of \cite{XiZa}).\\ 
\begin{lem}\label{XiZa2}
	\begin{enumerate}
		\item[(i)] If $q$ is fixed and $g\rightarrow\infty,$ then for $H\in \hh,$ the quantity \\
		$\log{N_{q}(J_{H})} -g\log{q}+\delta_{\gamma/2}\log(1- 1/q)$ converges weakly to a random variable $\R,$ whose characteristic function $\phi_{\R}(t)$ is given by 
		\begin{align*}
			\phi_{\R}(t)&=
			 1
			 + 
			 \sum\limits_{n=1}\limits^{\infty}
			 \frac{1}{2^{n}n!}
			 \sum\limits_{\substack{P_{1},...,P_{n}}}
			 \prod\limits_{j=1}\limits^{n}
			 \left(
			 \frac{(1-\mid P_{j} \mid^{-1})^{-it}
			 	+ (1+\mid P_{j} \mid^{-1})^{-it}-2}{(1+ \mid P_{j}\mid ^{-1})}
			 \right),   
		\end{align*}
		for all $t$ in $\mathbb{R}.$
		Here the inner sum is over distinct monic, irreducible polynomials $P_{1},...P_{n},$ in $\F_{q}[x].$  
		\item[(ii)] If both $q, g\rightarrow\infty,$ then over $\hh,$ the quantity
		$q^{1/2}\left(\log{N_{q}(J_{H})} -g\log{q}\right)$ has a standard Gaussian.
	\end{enumerate}
\end{lem}

We observe from \eqref{desingular} that, for $g \rightarrow\infty,$\\
\begin{equation*}
\log{N_{q}({\N})}-(4g-4)\log{q}=\log{N_{q^2}(J_{H})}-2g\log{q}.
\end{equation*}
In Lemma \eqref{XiZa2}, if we change the base field from $\F_{q}$ to $\F_{q^2}$, we get that the random variable
\[
\log{N_{q^2}(J_{H})} -2g\log{q}+\delta_{\gamma/2}\log(1- 1/q^2)
\]
converges weakly to a random variable $\R,$ whose characteristic function $\phi_{\R}(\tau)$ is as in Lemma \eqref{XiZa2}(i), such that, in the inner summation of the definition of $\phi_{\R}(t)$, distinct monic irreducible polynomials $P_{1},...,P_{n},$  are in $\F_{q^2}[x].$ Therefore, the random variable
\[
\log{N_{q}({\NN})}-(4g-4)\log{q} + \delta_{\gamma/2}\log(1- 1/q^2),
\]
 also converges weakly to $\R,$ and hence the proof of Theorem \eqref{thm1.6}$(1).$

		Next we see the case when both $q, g\rightarrow\infty.$ From Lemma \eqref{XiZa2}(ii), it follows that for $H:y^{2}=F(x)$ in $H_{\gamma, q^2},$ that is, for the the monic square-free polynomial $F(x)$ in $\F_{q^2}[x]$, the random variable
		\[
		q\left(\log{N_{q^2}(J_{H})} -2g\log{q}\right)
		\]
		 has a standard Gaussian distribution. Therefore we conclude that 
		 \[
		 q\left( \log{N_{q}({\NN})}-(4g-4)\log{q} \right)
		 \]
		  is distributed as a standard Gaussian over the family $\hh$, and hence the Theorem \eqref{thm1.6}$(2).$

\section{Distribution on $\b$}\label{distribution Higgs}

\subsection{Proof of Theorem \ref{thm1.1}:}
Let $X$ be a Galois curve of degree $N$ and genus $g\geq 2$ over $\F_{q}.$ We recall from \eqref{def:Ag2}, 
\begin{equation*}
A_{g,2}=A_{1}+A_{2}+A_{3},
\end{equation*}
where $A_{1}, A_{2}, A_{3}$ are as defined in Theorem \ref{Schiffman 1.6}.
From the definition of $A_{1}$ we can write 
\begin{equation*}
A_{1}=\frac{  q^{4g}     \prod\limits_{l=1}\limits^{2g} 
	(1-q^{-1/2}e(\theta_{l,X}))  (1-q^{-3/2}e(\theta_{l,X}))               }
{(q-1)(q^2-1)}.
\end{equation*}
Taking logarithm on both sides of the above equation we get
\begin{equation}\label{def:A}
\log(A_{1}) =4g\log {q}- \log{((q-1)^{2}(q+1))} +E_{1}+E_{2}
\end{equation}
where 
\begin{equation}\label{def:E1}
E_{1}=\sum\limits_{l=1}\limits^{2g}\log (1-q^{-3/2}e(\theta_{l,X})),
\end{equation}
and
\begin{equation}\label{def:E2}
E_{2}=\sum\limits_{l=1}\limits^{2g}\log (1-q^{-1/2}e(\theta_{l,X})).
\end{equation}

 Following a similar argument as in Lemma \ref{epsilon-Z5.1}, one can show the following results. For detailed proofs of these corresponding results one can also see \cite[Proposition 3.2]{ASA} and \cite[Section $3$]{XiZa} respectively. We have

\begin{equation}\label{log-zeta bound}
|E_{1}|
\le
(N-1)\left( \frac{1}{\sqrt{q}}+ \frac{1}{q}+ 
\frac{1}{q^{2}}\left\lbrace 4 
+\log{\left(\frac{2\log{\left( \frac{2g\sqrt{q}}{(1-q^{-3/2})(N-1)} \right) }}
	{3\log{q}}\right)}\right\rbrace\right),
\end{equation}
and
\begin{equation}\label{log-Jacobian bound}
|E_{2}|
\le
(N-1)
\left(
\log{\max\left\lbrace
	1,
	\frac{     \log{
			\left(\frac{7g}{(N-1)} \right)}}
	{   \log{q}   }
	\right\rbrace }
+3
\right).
\end{equation}


Next we analyse the term $A_{2}.$
Similarly as in the case $A_{1},$ we can write 
\begin{equation*}
A_{2}=\frac{  q^{2g}     \prod\limits_{l=1}\limits^{2g} 
	(1-q^{-1/2}e(\theta_{l,X}))  (1+q^{-1/2}e(\theta_{l,X}))               }
{4(q+1)}.
\end{equation*}
Taking logarithm on both sides of the above equation we get
\begin{equation}\label{def:B}
\log(A_{2}) =2g\log {q}- \log{(4(q+1))} + E_{2}+ E_{3}
\end{equation}
where, $E_{2}$ is defined as in \eqref{def:E2}
and
\begin{equation}\label{def:E3}
E_{3}= \sum\limits_{l=1}\limits^{2g}\log (1+q^{-1/2}e(\theta_{l,X})).
\end{equation}
It is easy to see that
\begin{equation}\label{third error}
|E_{3}|
\le
(N-1)
\left(
\log{\max\left\lbrace
	1,
	\frac{     \log{
			\left(\frac{7g}{(N-1)} \right)}}
	{   \log{q}   }
	\right\rbrace }
+3
\right).
\end{equation}

Next, we draw our attention to the term $A_{3},$ which we write as
\begin{equation*}
A_{3}=A_{3,1}+A_{3,2}+A_{3,3},
\end{equation*}
where
\begin{equation*}
A_{3,1}=\frac{   \prod\limits_{l=1}^{2g}(1-\alpha_{l})^{2}       }  
{     4(q-1)    },
\end{equation*}
and
\begin{equation*}
A_{3,2}=-\frac{   \prod\limits_{l=1}^{2g}(1-\alpha_{l})^{2}       }  
{     2(q-1)^{2}    },
\end{equation*}
and
\begin{equation*}
A_{3,3}=-\frac{   \prod\limits_{l=1}^{2g}(1-\alpha_{l})^{2}       }  
{     2(q-1)    }
\sum\limits_{l=1}^{2g}
\frac{1} {(1-\alpha_{l})   }
\end{equation*}

It is easy to check that the sum 
\begin{equation}
\left| 
\sum\limits_{l=1}^{2g}
\frac{1} {(1-\alpha_{l})   }
\right| 
=
\left| 
\sum\limits_{l=1}^{2g}
\frac{1} {(1-\sqrt{q} e(\theta_{l,X}))   }
\right| 
\leq
4g/\sqrt{q}.
\end{equation}

From \eqref{def:A}, and \eqref{def:B}, we observe that
\begin{equation*}
\log A_{2}-\log A_{1}=-2g\log q+
\log   \left( \frac{(q-1)^2}{4}  \right) 
+E_{3}-E_{1}.
\end{equation*}
Therefore, using \eqref{log-zeta bound}, and \eqref{third error}, we obtain
\begin{equation}\label{B/A}
\left| \frac{A_{2}}{A_{1}}\right| 
\le  
\left|
\exp(3(N-1))
\left(
\frac{\log{7g/(N-1)}}{\log{q}}
\right)^{N-1}q^{-2g+2}
\right|.
\end{equation}
Similarly, we have 
\begin{align*}
\frac{A_{3,1}}  {A_{1}}
&=
q^{-2g}
\frac{(q^2 -1)}{4}
\exp(E_{2}-E_{1}),\\
\frac{A_{3,2}}  {A_{1}}
&=
2q^{-2g}
(q+1)
\exp(E_{2}-E_{1}),\\
\frac{A_{3,3}}  {A_{1}}
&\leq
 2gq^{-2g-1/2}
(q^{2} - 1)
\exp(E_{2}-E_{1})
 \end{align*}

Therefore,

\begin{equation}\label{A_{3}/A_{1}}
\left|\frac{A_{3}}{A_{1}}\right|
=
O\left( 
gq^{-2g+3/2}\exp(3(N-1))
\left( \frac{\log g}{\log q}\right)^{N-1} 
\right) .
\end{equation}

\noindent
\textbf{Final steps of the proof of Theorem \ref{thm1.1}:} From Theorem \ref{Schiffman 1.1}, and Theorem \ref{Schiffman 1.6}, we have
\[
N_{q}(\b)=q^{(4g-3)}(A_{1}+A_{2}+A_{3}),
\]
Taking logarithm on both sides we have
\[
\log (N_{q}(\b))
=(4g-3)\log q
+\log A_{1}
+\log \left( 1+\frac{A_{2}+A_{3}}{A_{1}}\right) .
\]
Using \eqref{B/A}, and \eqref{A_{3}/A_{1}}, we observe that
 \begin{align*}
&\left|
\vv
-(4g-3)\log q
-\log A_{1}
\right| \\
&=
O 
\left( 
gq^{-2g+2}
\exp(3(N-1))
\left(
\frac{\log{7g/(N-1)}}{\log{q}}
\right)^{N-1}
\right). 
\end{align*}
Finally using \eqref{log-zeta bound}, and \eqref{log-Jacobian bound}, in \eqref{def:A}, the theorem follows upon simplifying.

\subsection{Proof of Theorem \ref{thm1.2}:}\label{distribution}

Next we focus our attention on the family of hyperelliptic curves $\H_{\gamma,q}.$ 
Suppose $d$ is any odd integer. Corresponding to a hyperelliptic curve $H,$ we use the notation $\a_{2,d}(H)$ to denote the moduli space of stable Higgs bundles of rank 2 and degree $d$ defined over $H.$

Proceeding as in the proof of Theorem \ref{thm1.1} we obtain
\begin{equation*}
\vw
-(4g-3)\log q
=\log A_{1}
+
\bigcirc
\left( 
gq^{-2g+2}
\frac{\log {g}}{\log{q}}
\right).
\end{equation*} 
Now using \eqref{e1} in \eqref{def:A}, for a fixed integer $Z,$ we can write
\begin{align*}
& \vw
 -(8g-6)\log{q}
+(1+\delta_{\gamma/2})\log{((1-1/q)^{2}(1+1/q))} \\
&= 
\bigtriangleup_{Z}(F)
+ \epsilon_{Z,F} 
+ \bigcirc
\left( 
gq^{-2g+2}
\frac{\log {g}}{\log{q}}
\right),
\end{align*}
where
\begin{equation}\label{mainterm}
\bigtriangleup_{Z}(F):= 
\sum\limits_{n\leq Z}
(q^{-n}+q^{-2n})
n^{-1}\sum\limits_{\substack{f \neq \infty  \\ \deg f=n}}
\Lambda(f)\left(\frac{F}{f}\right),
\end{equation}
and
\[ \epsilon_{Z, F} = 
-\sum\limits_{n>Z}
(
q^{-n/2}+q^{-3n/2}
)
n^{-1}
\sum\limits_{l=1}^{2g}
e(n\theta_{l,X})
-\delta_{\gamma/2}
\sum\limits_{n>Z}
\frac{q^{-n}+q^{-2n}}{n}.
\]
It is easy to see that 
\begin{equation}\label{upper-bounds}
\mid \bigtriangleup_{Z}(F)\mid
\leq 
\left( 
1+\frac{1}{q}
\right) 
\log{Z} \ 
\text{ and } \quad \mid\epsilon_{Z,F}\mid = \bigcirc\left(\frac{g}{Z} q^ {-Z/2}\right).
\end{equation}


{\bf The distribution function:}\\
Over the probability space $\hh,$ we define the random variable
\[
\mathcal{R}^{\a} : \H_{\gamma,q} \rightarrow \mathbb{R}
\]
such that
\[
\R^{\a}(H)=\left(\vw\right) -(8g-6)\log{q}.
\]
Using the definition \eqref{def:log constant}, we can write 
\begin{equation}\label{N_F}
\R^{\a}(H)
-C_{q}(2)+\delta_{\gamma/2}\log{\left( 1-1/q\right)} = \bigtriangleup_{Z}(F) + \epsilon_{Z}(F),
\end{equation}
where
\begin{equation}\label{error}
\epsilon_{Z}(F):= \epsilon_{Z,F} 
+
\bigcirc
\left( gq^{-2g+2}\frac{\log {g}}{\log{q}}
\right).
\end{equation}
Choose $Z= \left[ \frac{\gamma}{3}\right].$ If we put $k=0$ and $1$ in the definition \eqref{def: random variable sequence}, from \eqref{mainterm} we can identify the random variable
\begin{equation}\label{Higgs random variable}
\bigtriangleup_{Z}(F):= \R_{(\gamma, q)}^{(0)} + \R_{(\gamma, q)}^{(1)}. 
\end{equation}
Therefore, for some absolute constant $c>0,$ putting together the results from \eqref{upper-bounds},\eqref{error}, and\eqref{Higgs random variable} in \eqref{N_F}, we conclude
\[
\mathcal{R}^{\a} -C_{q}(2)+\delta_{\gamma/2}\log{\left(1-1/q \right)} = \R_{(\gamma, q)}^{(0)} + \R_{(\gamma, q)}^{(1)}
+
O
\left( q^{-cg}\right).
\]
Now, for a positive integer $n \leq \log \gamma,$ one can show that (cf. \cite[Theorem $3$]{XiZa}, and \cite[Theorem $1.2$]{ASA}) the $n$th moment of the random variables $\R_{(\gamma, q)}^{(k),}$ for $k=0,1$ is given by
\begin{equation*}\label{Delta-HH}
\E
\left( 
\left( 
\R_{(\gamma, q)}^{(k)}
\right) ^{n}
\right) 
=
H^{(k)}(n)+T 
\end{equation*}
where, 
\begin{equation*}
H^{(k)}(n) := 
\sum\limits_{\substack{m_i\geq 1\\ 1\leq i\leq n}} 
\prod\limits_{i=1}\limits^{n}q^{-(k+1)m_{i}}m_{i}^{-1}
\sum\limits_{\substack{\deg f_{i}=m_{i}\\ 1\leq i \leq n\\f_{1}f_{2}...f_{n}=h^{2}}}
\Lambda(f_{1})\Lambda(f_{2})....\Lambda(f_{n})
\prod\limits_{P\mid h} (1+ \mid P \mid ^{-1})^{-1},
\end{equation*}
and $T=O(q^{-c'\gamma})$ for some $c'>0.$
Therefore, if $q$ is fixed then for each fixed positive integer $n$,
\begin{equation*}
\lim\limits_{\gamma\rightarrow \infty} 
\E
\left( 
\left( 
\R_{(\gamma, q)}^{(k)} 
\right)^n 
\right) 
 = H^{(k)}(n).
\end{equation*}
 From Proposition\ref{H(r)1} we obtain the following result when $k=0$ and $1.$ For details one can also see \cite[Proposition $1$]{XiZa} for $k=0$, and \cite[Proposition $4.3$]{ASA} for $k=1.$

\begin{prop}\label{H(r)11}
	For any positive integer $n \geq 1$, we have for $k=0,1$
	\[
	H^{(k)}(n)= \sum_{s=1}^{n} \frac{n!}{2^{s}s!}
	\sum_{\substack{\sum\limits_{i=1}^{s}\lambda_{i}=n \\ \lambda_{i}\geq 1}} 
	\sum_{\substack{P_{1},...,P_{s}\\\text{distinct}}}
	\prod_{i=1}^{s}\frac{u_{P_{i}}^{\lambda_{i}}+ (-1)^{\lambda_{i}}v_{P_{i}}^{\lambda_{i}}}
	{\lambda_{i}!(1+ \mid P_{i}\mid ^{-1})}
	\]
	where the sum on the right hand side is over all positive integer 
	$\lambda_{i}, i= 1, 2, ..., s$ such that $\sum\limits_{i=1}\limits^{s}\lambda_{i}=n,$ 
	and over all distinct monic, irreducible polynomials $P_{i} \;\text{in}\; \F_{q}[x]$ with 
	\begin{eqnarray*}
		u_{P_{i}} &=& -\log (1-|P_{i}|^{-(k+1)}),\\
		v_{P_{i}} &=& \log(1+ |P_{i}|^{-(k+1)}).
	\end{eqnarray*}
\end{prop}

Next we proceed similarly as in the proof of Theorem \ref{thm1.4} to compute the limiting distribution function. 
We consider that $\R_{k}$ is a random variable such that for any positive integer $n,$
\[
\mathbb{E}\left( \R_{k}^n\right)=H^{(k)}(n).
\]
We know that characteristic function uniquely determines the distribution function. Therefore, for any real number $t,$ let the characteristic function $\phi_{\R_{k}}(t)$ of $\R_{k}$ is given by
\[
\phi_{\R_{k}}(t)
=
\mathbb{E}
\left(
e^{it\R_{k}}
\right). 
\]
Writing the characteristic function $\phi_{\R_k}(t)$ in terms of the $n^{\text{th}}$ moment $H^{(k)}(n),$ and applying Proposition \ref{H(r)11}, we obtain Theorem \ref{thm1.2}(1).

Using similar procedure as in Proposition \ref{thm covariance} one can show that for a fixed $q$ and $g \rightarrow \infty, $ Theorem \ref{thm1.2}(2) holds.\\

The next result gives an asymptotic formula for $H^{(k)}(n)$ when $q\rightarrow \infty,$ for $k=0$ and $1.$
\begin{prop}\label{H(r)3}
	As $q\rightarrow \infty,$ for a fixed positive integer $n,$ we obtain 
	\[
	H^{(k)}(n)
	= 
	\frac{\delta_{n/2}n!}{2^{n/2}(n/2)!}
	q^{\frac{-(2k+1)n}{2}} 
	+ \bigcirc_{n}
	\left( q^{ \frac{-(2k+1)(n+1)}{2}} 
	\right)
	\] 
	for $k=0,1.$
\end{prop}
\begin{proof}
	This is an exact analogue of \cite[Proposition $3$]{XiZa} for $k=0,$ and \cite[Proposition $4.4$]{ASA} for $k=1,$ and again we skip the proof. 
\end{proof}

For each $k=0$ and $1,$ considering 
$q^{(2k+1)/2}\R_{k}$ as a random variable on the space $\hh$, as both $\gamma,$ and $q \rightarrow \infty,$ we see that
all its moments are asymptotic to the corresponding moments of a 
standard Gaussian distribution where the odd moments vanish and the even moments are 
\begin{equation*}
\frac{1}{\sqrt{2\pi}}\int_{-\infty}^{\infty}\tau^{2n}e^{-\tau^{2}/2} d\tau = \frac{(2n)!}{2^{n}n!}.
\end{equation*} 
Hence the corresponding characteristic function converges to characteristic function of Gaussian distribution for each random variable $\R_{k}$. Finally using Continuity theorem ( see Theorem $3.3.6$ in \cite{Du}), we obtain result $(3)$ of Theorem \ref{thm1.2}.

\vspace{2cm}
\textbf{Acknowledgements:} We would like to thank D. S. Nagaraj for his valuable suggestions and comments when the results of this paper were presented to him. We also express our sincere gratitude to U.N. Bhosle for helping us in understanding the paper \cite{DeRa}, which was a great help while preparing this work. The second author also thanks Suhas B. N. for having many useful discussions, which helped her understand the various aspects of moduli spaces. We would like to thank O. Schiffmann for some useful email communication.


\begin{thebibliography}{10}
	
	
	
	\bibitem{B1}
	Balaji, V., \emph{Cohomology of certain moduli spaces of vector bundles}, Proc. Indian Acad. Sci. Math. Sci.  {\bf 98} (1988), 1--24.
	
	
	\bibitem{BaSe1}
	Balaji, V. and Seshadri, C. S., \emph{ Poincar\'{e} polynomials of some moduli varieties}, Algebraic geometry and analytic geometry ({T}okyo, 1990). Springer, Tokyo  (1991), 1--25.
	
	
	\bibitem{DeRa}
	Desale, U. V. and Ramanan, S., \emph{ Poincar\'{e} polynomials of the variety of stable bundles}, Math. Ann.  {\bf 216} (1975), 233--244.
	
	


\bibitem{Deu}
Deuring, M., \emph{ Die {T}ypen der {M}ultiplikatorenringe elliptischer
	{F}unktionenk\"{o}rper}, Abh. Math. Sem. Hansischen Univ.  {\bf 14} (1941), 197--272.


\bibitem{ASA}
Dey, A. and Dey, S. and Mukhopadhyay, A., \emph{ Statistics of moduli space of vector bundles}, Bull. Sci. Math.  {\bf 151} (2019), 13--33.


\bibitem{Du} Durrett, R., \emph{Probability theory and examples}, Fourth edition, Cambridge Series in Statistical and Probabilistic Mathematics, Cambridge University Press, Cambridge {\bf 31}(2010).  


\bibitem{FaRu}
Faifman, D. and Rudnick, Z., \emph{Statistics of the zeroes of zeta functions in families of hyperelliptic curves over a finite field}, Compositio Math.     {\bf 146} (2010), 81--101.




\bibitem{GhLe}
Ghione, F. and Letizia, M., \emph{ Effective divisors of higher rank on a curve and the {S}iegel
	formula}, Compositio Math.  {\bf 83} (1992), 147--159.


\bibitem{HaNa} Harder, G. and Narasimhan, M. S., \emph{ On the cohomology groups of moduli spaces of vector bundles on curves}, Math.ann.  {\bf 212} (1975), 215--248.




\bibitem{HRF08}
Hausel, T. and Rodriguez-Villegas, F., \emph{ Mixed {H}odge polynomials of character varieties}, Invent. Math.  {\bf 174} (2008), 555--624.


\bibitem{Hitchin}
Hitchin, N. J., \emph{ The self-duality equations on a {R}iemann surface}, Proc. London Math. Soc. (3).  {\bf 55} (1987), 59--126.


\bibitem{KaSa} N. M. Katz and P. Sarnak, \emph{Random matrices, Frobenius eigenvalues, and monodromy,} American Mathematical Society Colloquium Publications,  American Mathematical Society, Providence, RI, {\bf 45}(1999).



\bibitem{LaNgo}
Laumon, G. and Ng\^{o}, B. C., \emph{ Le lemme fondamental pour les groupes unitaires}, Ann. of Math. (2).  {\bf 168} (2008), 477--573.



\bibitem{Le}
Le Potier, J., \emph{ Lectures on vector bundles}, Cambridge Studies in Advanced Mathematics.  {\bf 54} Cambridge University Press, Cambridge (1997).





\bibitem{MeSe}
Mehta, V. B. and Seshadri, C. S., \emph{ Moduli of vector bundles on curves with parabolic structures}, Math. Ann.  {\bf 248} (1980), 205--239.




\bibitem{Mellit20}
Mellit, A., \emph{ Poincar\'{e} polynomials of character varieties, {M}acdonald
	polynomials and affine {S}pringer fibers}, Ann. of Math. (2).  {\bf 192} (2020), 165--228.



\bibitem{Moz12}
Mozgovoy, S., \emph{ Solutions of the motivic {ADHM} recursion formula}, Int. Math. Res. Not.(2012), no. 18,  4218--4244.




\bibitem{NaRa1}
Narasimhan, M. S. and Ramanan, S., \emph{ Moduli of vector bundles on a compact {R}iemann surface}, Ann. of Math. (2).  {\bf 89} (1969), 14--51.




\bibitem{NaSe}
Narasimhan, M. S. and Seshadri, C. S., \emph{ Stable and unitary vector bundles on a compact {R}iemann
	surface}, Ann. of Math. (2).  {\bf 82} (1965), 540--567.

\bibitem{Ne}
Newstead, P. E., \emph{Introduction to moduli problems and orbit spaces}, Tata Institute of Fundamental Research Lectures on Mathematics
and Physics, Tata Institute of Fundamental Research, Bombay; by the Narosa Publishing House, New Delhi {\bf 51}(1978).


\bibitem{Ngo}
Ng\^{o}, B. C., \emph{ Fibration de {H}itchin et endoscopie}, Invent. Math.  {\bf 164} (2006), 399--453.



\bibitem{Qu}
Quebbemann, H. G., \emph{ Estimates of regulators and class numbers in function fields}, J. Reine Angew. Math.  {\bf 419} (1991), 79--87.


\bibitem{Ro}
Rosen, M., \emph{Number theory in function fields, Graduate texts in Mathematics} Springer-Verlag, New York {\bf 210}(2002).



\bibitem{RoTs}
Rosenbloom, M. Y. and Tsfasman, M. A., \emph{ Multiplicative lattices in global fields}, Invent. Math.  {\bf 101} (1990), 687--696.


\bibitem{Ru}
Rudnick, Z., \emph{ Traces of high powers of the Frobenious class in the hyperelliptic ensemble}, Acta Arithmetica  {\bf 143} (2010), 81--99.


\bibitem{Sch}
Schiffmann, O., \emph{ Indecomposable vector bundles and stable {H}iggs bundles over
	smooth projective curves}, Ann. of Math. (2).  {\bf 183} (2016), 297--362.



\bibitem{Se2}
Seshadri, C. S., \emph{ Space of unitary vector bundles on a compact {R}iemann
		surface}, Ann. of Math. (2).  {\bf 85} (1967), 303--336.




\bibitem{Se4}
Seshadri, C. S., \emph{ Quotient spaces modulo reductive algebraic groups}, Ann. of Math. (2).  {\bf 95} (1972), 511--5562.


\bibitem{Se}
Seshadri, C. S., \emph{ Desingularisation of the moduli varieties of vector bundles on
	curves}, Proceedings of the {I}nternational {S}ymposium on {A}lgebraic
{G}eometry ({K}yoto {U}niv., {K}yoto, 1977),  (1978), 155--184.

\bibitem{Sh}
Shparlinski, I., \emph{ On the size of the {J}acobians of curves over finite fields}, Bull. Braz. Math. Soc. (N.S.).  {\bf 39} (2008), 587--595.



\bibitem{StTe}
Stein, A. and Teske, E., \emph{ Explicit bounds and heuristics on class numbers in
	hyperelliptic function fields}, Math. Comp.  {\bf 71} (2002), 837--861.



\bibitem{Ts}
Tsfasman, M. A., \emph{ Some remarks on the asymptotic number of points}, Coding theory and algebraic geometry ({L}uminy, 1991).  {\bf 1518} (1992), 178--192.





\bibitem{XiZa}
Xiong, M. and Zaharescu, A., \emph{ Statistics of the {J}acobians of hyperelliptic curves over finite fields}, Math. Res. Lett.  {\bf 19} (2012), 255--272.




\end{thebibliography}

	\end{document}